 
\hsize=13.50cm     
\vsize=18cm       
\parindent=12pt   \parskip=0pt      
\pageno=1 


\hoffset=15mm    
\voffset=1cm    
 

\ifnum\mag=\magstep1
\hoffset=-2mm   
\voffset=.8cm   
\fi


\pretolerance=500 \tolerance=1000  \brokenpenalty=5000

\catcode`\@=11

\font\eightrm=cmr8         \font\eighti=cmmi8
\font\eightsy=cmsy8        \font\eightbf=cmbx8
\font\eighttt=cmtt8        \font\eightit=cmti8
\font\eightsl=cmsl8        \font\sixrm=cmr6
\font\sixi=cmmi6           \font\sixsy=cmsy6
\font\sixbf=cmbx6


\font\tengoth=eufm10       \font\tenbboard=msbm10
\font\eightgoth=eufm10 at 8pt      \font\eightbboard=msbm10 at 8 pt
\font\sevengoth=eufm7      \font\sevenbboard=msbm7
\font\sixgoth=eufm7 at 6 pt        \font\fivegoth=eufm5

 \font\tencyr=wncyr10       
\font\eightcyr=wncyr10 at 8 pt      
\font\sevencyr=wncyr10 at 7 pt      
\font\sixcyr=wncyr10 at 6 pt


\skewchar\eighti='177 \skewchar\sixi='177
\skewchar\eightsy='60 \skewchar\sixsy='60


\newfam\gothfam           \newfam\bboardfam
\newfam\cyrfam

\def\tenpoint{%
  \textfont0=\tenrm \scriptfont0=\sevenrm \scriptscriptfont0=\fiverm
  \def\rm{\fam\z@\tenrm}%
  \textfont1=\teni  \scriptfont1=\seveni  \scriptscriptfont1=\fivei
  \def\oldstyle{\fam\@ne\teni}\let\old=\oldstyle
  \textfont2=\tensy \scriptfont2=\sevensy \scriptscriptfont2=\fivesy
  \textfont\gothfam=\tengoth \scriptfont\gothfam=\sevengoth
  \scriptscriptfont\gothfam=\fivegoth
  \def\goth{\fam\gothfam\tengoth}%
  \textfont\bboardfam=\tenbboard \scriptfont\bboardfam=\sevenbboard
  \scriptscriptfont\bboardfam=\sevenbboard
  \def\bb{\fam\bboardfam\tenbboard}%
 \textfont\cyrfam=\tencyr \scriptfont\cyrfam=\sevencyr
  \scriptscriptfont\cyrfam=\sixcyr
  \def\cyr{\fam\cyrfam\tencyr}%
  \textfont\itfam=\tenit
  \def\it{\fam\itfam\tenit}%
  \textfont\slfam=\tensl
  \def\sl{\fam\slfam\tensl}%
  \textfont\bffam=\tenbf \scriptfont\bffam=\sevenbf
  \scriptscriptfont\bffam=\fivebf
  \def\bf{\fam\bffam\tenbf}%
  \textfont\ttfam=\tentt
  \def\tt{\fam\ttfam\tentt}%
  \abovedisplayskip=12pt plus 3pt minus 9pt
  \belowdisplayskip=\abovedisplayskip
  \abovedisplayshortskip=0pt plus 3pt
  \belowdisplayshortskip=4pt plus 3pt 
  \smallskipamount=3pt plus 1pt minus 1pt
  \medskipamount=6pt plus 2pt minus 2pt
  \bigskipamount=12pt plus 4pt minus 4pt
  \normalbaselineskip=12pt
  \setbox\strutbox=\hbox{\vrule height8.5pt depth3.5pt width0pt}%
  \let\bigf@nt=\tenrm       \let\smallf@nt=\sevenrm
  \normalbaselines\rm}

\def\eightpoint{%
  \textfont0=\eightrm \scriptfont0=\sixrm \scriptscriptfont0=\fiverm
  \def\rm{\fam\z@\eightrm}%
  \textfont1=\eighti  \scriptfont1=\sixi  \scriptscriptfont1=\fivei
  \def\oldstyle{\fam\@ne\eighti}\let\old=\oldstyle
  \textfont2=\eightsy \scriptfont2=\sixsy \scriptscriptfont2=\fivesy
  \textfont\gothfam=\eightgoth \scriptfont\gothfam=\sixgoth
  \scriptscriptfont\gothfam=\fivegoth
  \def\goth{\fam\gothfam\eightgoth}%
  \textfont\cyrfam=\eightcyr \scriptfont\cyrfam=\sixcyr
  \scriptscriptfont\cyrfam=\sixcyr
  \def\cyr{\fam\cyrfam\eightcyr}%
  \textfont\bboardfam=\eightbboard \scriptfont\bboardfam=\sevenbboard
  \scriptscriptfont\bboardfam=\sevenbboard
  \def\bb{\fam\bboardfam}%
  \textfont\itfam=\eightit
  \def\it{\fam\itfam\eightit}%
  \textfont\slfam=\eightsl
  \def\sl{\fam\slfam\eightsl}%
  \textfont\bffam=\eightbf \scriptfont\bffam=\sixbf
  \scriptscriptfont\bffam=\fivebf
  \def\bf{\fam\bffam\eightbf}%
  \textfont\ttfam=\eighttt
  \def\tt{\fam\ttfam\eighttt}%
  \abovedisplayskip=9pt plus 3pt minus 9pt
  \belowdisplayskip=\abovedisplayskip
  \abovedisplayshortskip=0pt plus 3pt
  \belowdisplayshortskip=3pt plus 3pt 
  \smallskipamount=2pt plus 1pt minus 1pt
  \medskipamount=4pt plus 2pt minus 1pt
  \bigskipamount=9pt plus 3pt minus 3pt
  \normalbaselineskip=9pt
  \setbox\strutbox=\hbox{\vrule height7pt depth2pt width0pt}%
  \let\bigf@nt=\eightrm     \let\smallf@nt=\sixrm
  \normalbaselines\rm}

\tenpoint


\def\pc#1{\bigf@nt#1\smallf@nt}         \def\pd#1 {{\pc#1} }


\catcode`\;=\active
\def;{\relax\ifhmode\ifdim\lastskip>\z@\unskip\fi
\kern\fontdimen2  -1.2 \fontdimen3 \string;}

\catcode`\:=\active
\def:{\relax\ifhmode\ifdim\lastskip>\z@\unskip\fi\penalty\@M\ \fi\string:}

\catcode`\!=\active
\def!{\relax\ifhmode\ifdim\lastskip>\z@
\unskip\fi\kern\fontdimen2  -1.1 \fontdimen3 \string!}

\catcode`\?=\active
\def?{\relax\ifhmode\ifdim\lastskip>\z@
\unskip\fi\kern\fontdimen2  -1.1 \fontdimen3 \string?}

\def\^#1{\if#1i{\accent"5E\i}\else{\accent"5E #1}\fi}
\def\"#1{\if#1i{\accent"7F\i}\else{\accent"7F #1}\fi}

\frenchspacing


\newtoks\auteurcourant      \auteurcourant={\hfil}
\newtoks\titrecourant       \titrecourant={\hfil}

\newtoks\hautpagetitre      \hautpagetitre={\hfil}
\newtoks\baspagetitre       \baspagetitre={\hfil}

\newtoks\hautpagegauche     
\hautpagegauche={\eightpoint\rlap{\folio}\hfil\the\auteurcourant\hfil}
\newtoks\hautpagedroite     
\hautpagedroite={\eightpoint\hfil\the\titrecourant\hfil\llap{\folio}}

\newtoks\baspagegauche      \baspagegauche={\hfil} 
\newtoks\baspagedroite      \baspagedroite={\hfil}

\newif\ifpagetitre          \pagetitretrue  


\headline={\ifpagetitre\the\hautpagetitre
\else\ifodd\pageno\the\hautpagedroite\else\the\hautpagegauche\fi\fi}

\footline={\ifpagetitre\the\baspagetitre\else
\ifodd\pageno\the\baspagedroite\else\the\baspagegauche\fi\fi
\global\pagetitrefalse}


\def\raggedbottom{\topskip 10pt plus 36pt\r@ggedbottomtrue}



\def\pointir{\unskip . --- \ignorespaces}


\def\Bigbreak{\vskip-\lastskip\bigbreak}
\def\Medbreak{\vskip-\lastskip\medbreak}


\def\ctexte#1\endctexte{%
  \hbox{$\vcenter{\halign{\hfill##\hfill\crcr#1\crcr}}$}}


\long\def\ctitre#1\endctitre{%
    \ifdim\lastskip<24pt\vskip-\lastskip\bigbreak\bigbreak\fi
  		\vbox{\parindent=0pt\leftskip=0pt plus 1fill
          \rightskip=\leftskip
          \parfillskip=0pt\bf#1\par}
    \bigskip\nobreak}

\long\def\section#1\endsection{%
\vskip 0pt plus 3\normalbaselineskip
\penalty-250
\vskip 0pt plus -3\normalbaselineskip
\Bigbreak
\message{[section \string: #1]}{\bf#1\unskip}\pointir}

\long\def\sectiona#1\endsection{%
\vskip 0pt plus 3\normalbaselineskip
\penalty-250
\vskip 0pt plus -3\normalbaselineskip
\Bigbreak
\message{[sectiona \string: #1]}%
{\bf#1}\medskip\nobreak}

\long\def\subsection#1\endsubsection{%
\Medbreak
{\it#1\unskip}\pointir}

\long\def\subsectiona#1\endsubsection{%
\Medbreak
{\it#1}\par\nobreak}

\def\rem#1\endrem{%
\Medbreak
{\it#1\unskip} : }

\def\remp#1\endrem{%
\Medbreak
{\pc #1\unskip}\pointir}

\def\rema#1\endrem{%
\Medbreak
{\it #1}\par\nobreak}

\def\newparwithcolon#1\endnewparwithcolon{
\Medbreak
{#1\unskip} : }

\def\newparwithpointir#1\endnewparwithpointir{
\Medbreak
{#1\unskip}\pointir}

\def\newpara#1\endnewpar{
\Medbreak
{#1\unskip}\smallskip\nobreak}


\long\def\th#1 #2\enonce#3\endth{%
   \Medbreak
   {\pc#1} {#2\unskip}\pointir{\it #3}\medskip}

\long\def\tha#1 #2\enonce#3\endth{%
   \Medbreak
   {\pc#1} {#2\unskip}\par\nobreak{\it #3}\medskip}


\long\def\Th#1 #2 #3\enonce#4\endth{%
   \Medbreak
   #1 {\pc#2} {#3\unskip}\pointir{\it #4}\medskip}

\long\def\Tha#1 #2 #3\enonce#4\endth{%
   \Medbreak
   #1 {\pc#2} #3\par\nobreak{\it #4}\medskip}


\def\decale#1{\smallbreak\hskip 28pt\llap{#1}\kern 5pt}
\def\decaledecale#1{\smallbreak\hskip 34pt\llap{#1}\kern 5pt}
\def\puce{\smallbreak\hskip 6pt{$\scriptstyle\bullet$}\kern 5pt}



\def\displaylinesno#1{\displ@y\halign{
\hbox to\displaywidth{$\@lign\hfil\displaystyle##\hfil$}&
\llap{$##$}\crcr#1\crcr}}


\def\ldisplaylinesno#1{\displ@y\halign{ 
\hbox to\displaywidth{$\@lign\hfil\displaystyle##\hfil$}&
\kern-\displaywidth\rlap{$##$}\tabskip\displaywidth\crcr#1\crcr}}


\def\eqalign#1{\null\,\vcenter{\openup\jot\m@th\ialign{
\strut\hfil$\displaystyle{##}$&$\displaystyle{{}##}$\hfil
&&\quad\strut\hfil$\displaystyle{##}$&$\displaystyle{{}##}$\hfil
\crcr#1\crcr}}\,}


\def\system#1{\left\{\null\,\vcenter{\openup1\jot\m@th
\ialign{\strut$##$&\hfil$##$&$##$\hfil&&
        \enskip$##$\enskip&\hfil$##$&$##$\hfil\crcr#1\crcr}}\right.}


\let\@ldmessage=\message

\def\message#1{{\def\pc{\string\pc\space}%
                \def\'{\string'}\def\`{\string`}%
                \def\^{\string^}\def\"{\string"}%
                \@ldmessage{#1}}}



\def\up#1{\raise 1ex\hbox{\smallf@nt#1}}


\def\qed{\raise -2pt\hbox{\vrule\vbox to 10pt{\hrule width 4pt
                 \vfill\hrule}\vrule}}

\def\cqfd{\unskip\penalty 500\quad\qed\medbreak}

\def\virg{\raise .4ex\hbox{,}}   


\def\build#1_#2^#3{\mathrel{
\mathop{\kern 0pt#1}\limits_{#2}^{#3}}}


\def\boxit#1#2{%
\setbox1=\hbox{\kern#1{#2}\kern#1}%
\dimen1=\ht1 \advance\dimen1 by #1 \dimen2=\dp1 \advance\dimen2 by #1 
\setbox1=\hbox{\vrule height\dimen1 depth\dimen2\box1\vrule}%
\setbox1=\vbox{\hrule\box1\hrule}%
\advance\dimen1 by .6pt \ht1=\dimen1 
\advance\dimen2 by .6pt \dp1=\dimen2  \box1\relax}


\catcode`\@=12

\showboxbreadth=-1  \showboxdepth=-1



    \input amssym.def
\input amssym.tex 

\def\Grille{\setbox13=\vbox to 5\unitlength{\hrule width 109mm\vfill} 
\setbox13=\vbox to 65mm{\offinterlineskip\leaders\copy13\vfill\kern-1pt\hrule} 
\setbox14=\hbox to 5\unitlength{\vrule height 65mm\hfill} 
\setbox14=\hbox to 109mm{\leaders\copy14\hfill\kern-2mm\vrule height 65mm}
\ht14=0pt\dp14=0pt\wd14=0pt \setbox13=\vbox to
0pt{\vss\box13\offinterlineskip\box14} \wd13=0pt\box13}


\def\fleche(#1,#2)\dir(#3,#4)\long#5{%
\noalign{\leftput(#1,#2){\vector(#3,#4){#5}}}}

\def\ligne(#1,#2)\dir(#3,#4)\long#5{%
\noalign{\leftput(#1,#2){\lline(#3,#4){#5}}}}

\def\put(#1,#2)#3{\noalign{\setbox1=\hbox{%
    \kern #1\unitlength
    \raise #2\unitlength\hbox{$#3$}}%
    \ht1=0pt \wd1=0pt \dp1=0pt\box1}}


\def\diagram#1{\def\normalbaselines{\baselineskip=0pt\lineskip=5pt}
\matrix{#1}}

\def\hfl#1#2#3{\smash{\mathop{\hbox to#3{\rightarrowfill}}\limits
^{\scriptstyle#1}_{\scriptstyle#2}}}

\def\gfl#1#2#3{\smash{\mathop{\hbox to#3{\leftarrowfill}}\limits
^{\scriptstyle#1}_{\scriptstyle#2}}}

\def\vfl#1#2#3{\llap{$\scriptstyle #1$}
\left\downarrow\vbox to#3{}\right.\rlap{$\scriptstyle #2$}}

\def\ufl#1#2#3{\llap{$\scriptstyle #1$}
\left\uparrow\vbox to#3{}\right.\rlap{$\scriptstyle #2$}}

\def\pafl#1#2#3{\llap{$\scriptstyle #1$}
\left\Vert\vbox to#3{}\right.\rlap{$\scriptstyle #2$}}


 \message{`lline' & `vector' macros from LaTeX}
 \catcode`@=11
\def\{{\relax\ifmmode\lbrace\else$\lbrace$\fi}
\def\}{\relax\ifmmode\rbrace\else$\rbrace$\fi}
\def\newcount{\alloc@0\count\countdef\insc@unt}
\def\newdimen{\alloc@1\dimen\dimendef\insc@unt}
\def\newwrite{\alloc@7\write\chardef\sixt@@n}

\newwrite\@unused
\newcount\@tempcnta
\newcount\@tempcntb
\newdimen\@tempdima
\newdimen\@tempdimb
\newbox\@tempboxa

\def\@spaces{\space\space\space\space}
\def\@whilenoop#1{}
\def\@whiledim#1\do #2{\ifdim #1\relax#2\@iwhiledim{#1\relax#2}\fi}
\def\@iwhiledim#1{\ifdim #1\let\@nextwhile=\@iwhiledim
        \else\let\@nextwhile=\@whilenoop\fi\@nextwhile{#1}}
\def\@badlinearg{\@latexerr{Bad \string\line\space or \string\vector
   \space argument}}
\def\@latexerr#1#2{\begingroup
\edef\@tempc{#2}\expandafter\errhelp\expandafter{\@tempc}%
\def\@eha{Your command was ignored.
^^JType \space I <command> <return> \space to replace it
  with another command,^^Jor \space <return> \space to continue without
it.} 
\def\@ehb{You've lost some text. \space \@ehc}
\def\@ehc{Try typing \space <return>
  \space to proceed.^^JIf that doesn't work, type \space X <return> \space to
  quit.}
\def\@ehd{You're in trouble here.  \space\@ehc}

\typeout{LaTeX error. \space See LaTeX manual for explanation.^^J
 \space\@spaces\@spaces\@spaces Type \space H <return> \space for
 immediate help.}\errmessage{#1}\endgroup}
\def\typeout#1{{\let\protect\string\immediate\write\@unused{#1}}} 

\font\tenln    = line10
\font\tenlnw   = linew10

\newdimen\@wholewidth
\newdimen\@halfwidth
\newdimen\unitlength 

\unitlength =1pt


\def\thinlines{\let\@linefnt\tenln \let\@circlefnt\tencirc
  \@wholewidth\fontdimen8\tenln \@halfwidth .5\@wholewidth}
\def\thicklines{\let\@linefnt\tenlnw \let\@circlefnt\tencircw
  \@wholewidth\fontdimen8\tenlnw \@halfwidth .5\@wholewidth}

\def\linethickness#1{\@wholewidth #1\relax \@halfwidth .5\@wholewidth}

\newif\if@negarg

\def\lline(#1,#2)#3{\@xarg #1\relax \@yarg #2\relax
\@linelen=#3\unitlength
\ifnum\@xarg =0 \@vline
  \else \ifnum\@yarg =0 \@hline \else \@sline\fi
\fi}

\def\@sline{\ifnum\@xarg< 0 \@negargtrue \@xarg -\@xarg \@yyarg -\@yarg
  \else \@negargfalse \@yyarg \@yarg \fi
\ifnum \@yyarg >0 \@tempcnta\@yyarg \else \@tempcnta -\@yyarg \fi
\ifnum\@tempcnta>6 \@badlinearg\@tempcnta0 \fi
\setbox\@linechar\hbox{\@linefnt\@getlinechar(\@xarg,\@yyarg)}%
\ifnum \@yarg >0 \let\@upordown\raise \@clnht\z@
   \else\let\@upordown\lower \@clnht \ht\@linechar\fi
\@clnwd=\wd\@linechar
\if@negarg \hskip -\wd\@linechar \def\@tempa{\hskip -2\wd\@linechar}\else
     \let\@tempa\relax \fi
\@whiledim \@clnwd <\@linelen \do
  {\@upordown\@clnht\copy\@linechar
   \@tempa
   \advance\@clnht \ht\@linechar
   \advance\@clnwd \wd\@linechar}%
\advance\@clnht -\ht\@linechar
\advance\@clnwd -\wd\@linechar
\@tempdima\@linelen\advance\@tempdima -\@clnwd
\@tempdimb\@tempdima\advance\@tempdimb -\wd\@linechar
\if@negarg \hskip -\@tempdimb \else \hskip \@tempdimb \fi
\multiply\@tempdima \@m
\@tempcnta \@tempdima \@tempdima \wd\@linechar \divide\@tempcnta \@tempdima
\@tempdima \ht\@linechar \multiply\@tempdima \@tempcnta
\divide\@tempdima \@m
\advance\@clnht \@tempdima
\ifdim \@linelen <\wd\@linechar
   \hskip \wd\@linechar
  \else\@upordown\@clnht\copy\@linechar\fi}

\def\@hline{\ifnum \@xarg <0 \hskip -\@linelen \fi
\vrule height \@halfwidth depth \@halfwidth width \@linelen
\ifnum \@xarg <0 \hskip -\@linelen \fi}

\def\@getlinechar(#1,#2){\@tempcnta#1\relax\multiply\@tempcnta 8
\advance\@tempcnta -9 \ifnum #2>0 \advance\@tempcnta #2\relax\else
\advance\@tempcnta -#2\relax\advance\@tempcnta 64 \fi
\char\@tempcnta}

\def\vector(#1,#2)#3{\@xarg #1\relax \@yarg #2\relax
\@linelen=#3\unitlength
\ifnum\@xarg =0 \@vvector
  \else \ifnum\@yarg =0 \@hvector \else \@svector\fi
\fi} 

\def\@hvector{\@hline\hbox to 0pt{\@linefnt
\ifnum \@xarg <0 \@getlarrow(1,0)\hss\else
    \hss\@getrarrow(1,0)\fi}}

\def\@vvector{\ifnum \@yarg <0 \@downvector \else \@upvector \fi}

\def\@svector{\@sline
\@tempcnta\@yarg \ifnum\@tempcnta <0 \@tempcnta=-\@tempcnta\fi
\ifnum\@tempcnta <5
  \hskip -\wd\@linechar
  \@upordown\@clnht \hbox{\@linefnt  \if@negarg
  \@getlarrow(\@xarg,\@yyarg) \else \@getrarrow(\@xarg,\@yyarg) \fi}%
\else\@badlinearg\fi}

\def\@getlarrow(#1,#2){\ifnum #2 =\z@ \@tempcnta='33\else
\@tempcnta=#1\relax\multiply\@tempcnta \sixt@@n \advance\@tempcnta
-9 \@tempcntb=#2\relax\multiply\@tempcntb \tw@
\ifnum \@tempcntb >0 \advance\@tempcnta \@tempcntb\relax
\else\advance\@tempcnta -\@tempcntb\advance\@tempcnta 64
\fi\fi\char\@tempcnta}

\def\@getrarrow(#1,#2){\@tempcntb=#2\relax
\ifnum\@tempcntb < 0 \@tempcntb=-\@tempcntb\relax\fi
\ifcase \@tempcntb\relax \@tempcnta='55 \or
\ifnum #1<3 \@tempcnta=#1\relax\multiply\@tempcnta
24 \advance\@tempcnta -6 \else \ifnum #1=3 \@tempcnta=49
\else\@tempcnta=58 \fi\fi\or
\ifnum #1<3 \@tempcnta=#1\relax\multiply\@tempcnta
24 \advance\@tempcnta -3 \else \@tempcnta=51\fi\or
\@tempcnta=#1\relax\multiply\@tempcnta
\sixt@@n \advance\@tempcnta -\tw@ \else
\@tempcnta=#1\relax\multiply\@tempcnta
\sixt@@n \advance\@tempcnta 7 \fi\ifnum #2<0 \advance\@tempcnta 64 \fi
\char\@tempcnta}

\def\@vline{\ifnum \@yarg <0 \@downline \else \@upline\fi}

\def\@upline{\hbox to \z@{\hskip -\@halfwidth \vrule
  width \@wholewidth height \@linelen depth \z@\hss}}

\def\@downline{\hbox to \z@{\hskip -\@halfwidth \vrule
  width \@wholewidth height \z@ depth \@linelen \hss}}

\def\@upvector{\@upline\setbox\@tempboxa\hbox{\@linefnt\char'66}\raise
     \@linelen \hbox to\z@{\lower \ht\@tempboxa\box\@tempboxa\hss}}

\def\@downvector{\@downline\lower \@linelen
      \hbox to \z@{\@linefnt\char'77\hss}}

\thinlines

\newcount\@xarg
\newcount\@yarg
\newcount\@yyarg
\newcount\@multicnt
\newdimen\@xdim
\newdimen\@ydim
\newbox\@linechar
\newdimen\@linelen
\newdimen\@clnwd
\newdimen\@clnht
\newdimen\@dashdim
\newbox\@dashbox
\newcount\@dashcnt
 \catcode`@=12


\newbox\tbox
\newbox\tboxa

\def\leftzer#1{\setbox\tbox=\hbox to 0pt{#1\hss}%
     \ht\tbox=0pt \dp\tbox=0pt \box\tbox}

\def\rightzer#1{\setbox\tbox=\hbox to 0pt{\hss #1}%
     \ht\tbox=0pt \dp\tbox=0pt \box\tbox}

\def\centerzer#1{\setbox\tbox=\hbox to 0pt{\hss #1\hss}%
     \ht\tbox=0pt \dp\tbox=0pt \box\tbox}

%
\def\image(#1,#2)#3{\vbox to #1{\offinterlineskip
    \vss #3 \vskip #2}}


\def\leftput(#1,#2)#3{\setbox\tboxa=\hbox{%
    \kern #1\unitlength
    \raise #2\unitlength\hbox{\leftzer{#3}}}%
    \ht\tboxa=0pt \wd\tboxa=0pt \dp\tboxa=0pt\box\tboxa}

\def\rightput(#1,#2)#3{\setbox\tboxa=\hbox{%
    \kern #1\unitlength
    \raise #2\unitlength\hbox{\rightzer{#3}}}%
    \ht\tboxa=0pt \wd\tboxa=0pt \dp\tboxa=0pt\box\tboxa}

\def\centerput(#1,#2)#3{\setbox\tboxa=\hbox{%
    \kern #1\unitlength
    \raise #2\unitlength\hbox{\centerzer{#3}}}%
    \ht\tboxa=0pt \wd\tboxa=0pt \dp\tboxa=0pt\box\tboxa}

\unitlength=1mm

\magnification=\magstep1
\hsize=17,5truecm
\vsize=25.5truecm
\hoffset=-0.9truecm
\voffset=-0.8truecm
\topskip=1truecm
\footline={\tenrm\hfil\folio\hfil}
\raggedbottom
\abovedisplayskip=3mm 
\belowdisplayskip=3mm
\abovedisplayshortskip=0mm
\belowdisplayshortskip=2mm
\normalbaselineskip=12pt  
\normalbaselines

\def\Z{{\Bbb Z}}

\def\G{{\Bbb G}}

\def\br{{\rm Br}}
\def\k{{\overline k}}
\def\pic{{\rm Pic}}
\def\dim{{\rm dim}}
\def\gal{{\rm Gal}}
\def\hom{{\rm Hom}}

\def\X{{\cyr X}}
 \def\et{\hbox{\sevenrm \'et}}

{\it Groupe de Picard et groupe de Brauer des compactifications lisses d'espaces homog\`enes}

\bigskip


\medskip
Jean-Louis Colliot-Th\'el\`ene et Boris \`E. Kunyavski\u{\i}

\bigskip
\bigskip

{\bf Introduction}

\medskip Soient $k$ un corps de caract\'eristique nulle, ${\overline k}$ une cl\^oture
alg\'ebrique. Soit  $g$ le groupe de Galois de
${\overline k}$ sur $k$. Etant donn\'ee une $k$-vari\'et\'e alg\'ebrique $Z$, on note
${\overline Z}=Z \times_k{\k}$.  \`A toute  $k$-vari\'et\'e projective, lisse et
g\'eom\'etriquement connexe $Z$, 
on associe le $g$-module (module galoisien) d\'efini par le groupe de Picard ${\rm Pic}({\overline Z} )$.
Si la vari\'et\'e ${\overline Z}$ est unirationnelle, le groupe ${\rm Pic}({\overline Z} )$ est un groupe ab\'elien libre de type fini.
Si $Z_1$ et
$Z_2$ sont deux  $k$-vari\'et\'es projectives, lisses et
g\'eom\'etriquement connexes
$k$-birationnellement \'equivalentes, les modules galoisiens ${\rm Pic}({\overline Z}_1 )$ et
${\rm Pic}({\overline Z}_2)$ sont isomorphes \`a addition pr\`es de $g$-modules de permutation
(de type fini) ([CT/San2, Prop.~2.A.1]). Le groupe de cohomologie $H^1(g,{\rm Pic}({\overline Z} ))$
est un invariant $k$-birationnel \'etroitement li\'e au
 groupe $H^2_{\et}(Z,\G_m)$, qui est le groupe
de Brauer cohomologique de $Z$, not\'e $\br(Z)$ (voir la proposition 1.1 ci-dessous).

Voskresenski\u{\i} (1974)   montra que pour tout $k$-tore $T$ et toute $k$-compactification lisse $T_c$ de $T$, le module galoisien ${\rm Pic}({\overline T}_c )$ poss\`ede la propri\'et\'e remarquable  suivante : pour tout sous-groupe ferm\'e $h \subset g$, on a
${\rm Ext}^1_h({\rm Pic}({\overline T}_c),\Z)=0$ ([Vos1]; [Vos2], IV.4.49; voir aussi [CT/San1], \S 2 et [Vos3], \S 4.6).
Dans la terminologie introduite dans [CT/San1], le module galoisien  ${\rm Pic}({\overline T}_c )$
est un module flasque. 
Un  argument simple donn\'e dans un article r\'ecent ([B/K2], 2004) montre que ce r\'esultat
s'\'etend aux groupes lin\'eaires connexes : pour toute $k$-compactification lisse $G_c$ d'un $k$-groupe lin\'eaire connexe
$G$, le module galoisien ${\rm Pic}({\overline G}_c)$ est flasque.

Le fait que le module galoisien ${\rm Pic}({\overline G}_c)$ est flasque permet
de calculer ce module, \`a addition pr\`es d'un
$g$-module de permutation,  et de donner  une formule simple  pour le groupe de Brauer de $G_c$,
ceci  purement en termes du groupe $G$, sans avoir \`a en exhiber une compactification lisse  
  ([CT/San1] 
 dans le cas des tores;  [B/K2], [CT] pour les groupes lin\'eaires connexes).

\medskip 

{\it Dans cet article, nous montrons que des r\'esultats analogues valent plus g\'en\'eralement
pour les compactifications 
 lisses d'espaces homog\`enes de groupes lin\'eaires connexes,
 quand le stabilisateur g\'eom\'etrique
 est connexe.}

\medskip

Soient
$G$ un $k$-groupe 
connexe et $X/k$ un espace homog\`ene de $G$. Le {\it
stabilisateur g\'eom\'etrique}, c'est-\`a-dire le groupe d'isotropie d'un ${\overline k}$-point
de
${\overline X}=X \times_k {\overline k}$ est bien d\'efini \`a
$\k$-isomorphisme non unique pr\`es. On {\it note} ${\overline H}$ ce groupe. 
Supposons le groupe ${\overline H}$ connexe.
Il y a alors
un $k$-tore $T$ naturellement associ\'e au $G$-espace homog\`ene $X$,
tel que ${\overline T}$ soit le plus grand quotient torique ${\overline H}^{\rm tor}$ de ${\overline H}$.
Une d\'efinition de ce $k$-{\it tore associ\'e} est donn\'ee par Borovoi dans [Bo2, 4.1]. 
On en trouvera une autre au \S 1.

Soit $X_c$ une $k$-compactification lisse de $X$. La
$\k$-vari\'et\'e ${\overline X}_c$ est unirationnelle, le groupe de Picard ${\rm Pic}({\overline
X}_c)$ est  un $g$-module continu discret
$\Z$-libre de type fini et le groupe $\br({\overline X}_c)$ est fini.
On note $\br_1(X_c)$ le noyau de l'application
de restriction $\br(X_c) \to \br({\overline X}_c)$.
Le quotient du groupe de Brauer
$\br_1(X_c)$ par l'image  du groupe $\br(k)$ est un sous-groupe du groupe fini $H^1(g,{\rm
Pic}({\overline X}_c))$ (Prop. 1.1).

\medskip \`A tout $g$-module continu discret $M$ et tout entier naturel $i$ on associe le groupe
$$\X^{i}_{\omega}(k,M) = {\rm Ker} [H^{i}(g,M) \to  \prod_h H^{i}(h,M)],$$ o\`u $h$ parcourt les
sous-groupes ferm\'es procycliques de $g$.

\bigskip

Le but principal de l'article est d'\'etablir le th\'eor\`eme suivant (Th\'eor\`eme 5.1).

\medskip

{\bf Th\'eor\`eme A} {\it Soient $k$ un corps de caract\'eristique nulle,
$G$ un $k$-groupe lin\'eaire connexe,
$X$ une $k$-vari\'et\'e espace homog\`ene de
$G$, de stabilisateur g\'eom\'etrique connexe.  Soit $X_c$ une $k$-compactification lisse de $X$.

{\rm (i)}  Le $g$-module
${\rm Pic}({\overline X}_c)$ est un $g$-module flasque, c'est-\`a-dire que pour tout sous-groupe
ferm\'e $h \subset g$, on a $H^1(h,\hom_\Z({\rm Pic}({\overline X}_c),\Z))=0$, soit encore
${\rm Ext}^1_h({\rm Pic}({\overline X}_c),\Z)=0$.

{\rm (ii)}  Pour tout sous-groupe ferm\'e procyclique $h \subset g$, on a  $H^1(h,{\rm Pic}({\overline
X}_c))=0$.

{\rm (iii)  }  Soit 
$T$ le $k$-tore associ\'e au $G$-espace homog\`ene $X$, et soit $\hat T$ son groupe des caract\`eres. Si $G$ est un 
groupe lin\'eaire quasitrivial, i.e. extension d'un $k$-tore
quasitrivial par un $k$-groupe simplement connexe, alors le quotient 
 du groupe
$\br_1(X_c)$ par l'image du groupe $\br(k)$ s'injecte dans le groupe
$\X^{1}_{\omega}(k,{\hat T})$, et est isomorphe \`a ce dernier groupe si $X(k)\neq \emptyset$ ou
si $k$ est un corps de nombres.}

\medskip

 Sous l'hypoth\`ese de (iii), nous montrons comment 
  le $g$-module $\Z$-libre de type fini ${\rm Pic}({\overline X}_c)$ est
 d\'etermin\'e,  \`a addition pr\`es d'un
$g$-module de permutation, par le $k$-tore $T$ --   en particulier il ne d\'epend pas du
groupe quasi-trivial $G$.  
Ce th\'eor\`eme est une extension naturelle de r\'esultats 
 bien connus dans le cas o\`u $G$ est un tore alg\'ebrique  (\S 2).

Pour $G$ semi-simple simplement connexe, l'\'enonc\'e (iii) combin\'e avec un th\'eor\`eme de
Bogomolov (Th\'eor\`eme 1.4 ci-dessous) \'etablit  
une conjecture  propos\'ee  au \S 5 de [CT/K]. 
Cette conjecture \'etait motiv\'ee par les r\'esultats de 
Borovoi [Bo2] sur le principe de Hasse et l'approximation faible pour les espaces homog\`enes
d\'efinis sur un corps de nombres.

\bigskip

Un ingr\'edient important de la d\'emonstration est le th\'eor\`eme suivant (Th\'eor\`eme 4.2).

\smallskip
{\bf Th\'eor\`eme B} {\it Soit $A$ un anneau de valuation discr\`ete de corps des fractions $K$,
de corps r\'esiduel $k$ de caract\'eristique nulle. Soit $G$ un
$K$-groupe quasitrivial et soit $E/K$ un $G$-espace homog\`ene  de
stabilisateur g\'eom\'etrique connexe et de tore associ\'e trivial. Soit
$X$ un
$A$-sch\'ema propre, r\'egulier, int\`egre, dont la fibre g\'en\'erique contient $E$ comme
ouvert dense. Alors il existe une composante de multiplicit\'e 1 de la fibre sp\'eciale de $X/A$
qui est g\'eom\'etriquement int\`egre sur son corps de base $k$.}

\smallskip

L'hypoth\`ese que le tore $T$ associ\'e est trivial  ($T=1$)
 \'equivaut au fait que le quotient du
stabilisateur g\'eom\'etrique ${\overline H}$ par son radical unipotent est un groupe
semi-simple (connexe).

\medskip

Nous discutons aussi un autre invariant, \`a savoir l'ensemble des classes de
$R$-\'equivalence sur les $k$-points de $X_c$. C'est l'objet du \S 6,
o\`u l'on obtient une minoration de l'ensemble $X_c(k)/R$ lorsque le corps
$k$ est un ``bon'' corps de dimension cohomologique 2 (par exemple un corps $p$-adique).

\medskip

Dans une premi\`ere mouture du pr\'esent article, sous l'hypoth\`ese que
$G$ est semi-simple simplement connexe, nous avions \'etabli 
les points (ii) et (iii) du Th\'eor\`eme A,  et prouv\'e des cas particuliers du 
th\'eor\`eme B,
dont nous avions  montr\'e qu'il implique  le th\'eor\`eme A. Nos preuves  passaient  par un long
d\'etour arithm\'etique (r\'eduction au cas des corps de nombres et de leurs compl\'et\'es).
Dans cette version primitive de l'article, un outil de base
\'etait un th\'eor\`eme de Borovoi [Bo1] : si $k$ est un corps $p$-adique,
$G$ un $k$-groupe semi-simple simplement connexe et $X$ un espace homog\`ene de
$G$ de stabilisateur g\'eom\'etrique connexe sans quotient torique, alors
$X(k) \neq \emptyset$.
Lors d'une discussion o\`u nous mentionnions l'\'enonc\'e du Th\'eor\`eme B,
Ofer Gabber nous a sugg\'er\'e d'utiliser une variante purement alg\'ebrique du th\'eor\`eme de Borovoi
pour \'etablir le Th\'eor\`eme B, en utilisant une r\'eduction toute diff\'erente
de celle mentionn\'ee ci-dessus.

Nous remercions vivement O.~Gabber pour son importante suggestion, et
P.~Gille et M.~Borovoi pour diverses remarques. M. Borovoi a en
particulier ind\'ependamment  remarqu\'e que l'hypoth\`ese  $G$
semi-simple simplement connexe, que nous avions impos\'ee dans la version pr\'ec\'edente, 
peut \^etre omise de la plupart de nos
\'enonc\'es.

\bigskip

{\sevenrm Ce travail a \'et\'e commenc\'e lors d'un s\'ejour de 
B. Kunyavski{\u\i}  \`a l'Universit\'e  Paris-Sud en septembre 2004.
Les r\'esultats finaux ont \'et\'e expos\'es  \`a la conf\'erence  ``Applications of torsors to Galois cohomology and Lie theory'' qui s'est tenue au BIRS, \`a Banff (Alberta, Canada)
en avril 2005.
Pour sa recherche, B. Kunyavski\u\i  \ a aussi b\'en\'efici\'e du soutien partiel du
 Minist\`ere de l'Absorption (Isra\"el), de l'Institut Emmy Noether (Fondation Minerva), de la
Fondation Isra\'elienne des Sciences (ISF) fond\'ee par l'Acad\'emie Isra\'elienne des Sciences
et des Lettres (le programme Centre d'Excellence ``Group-theoretic Methods in the Study of
Algebraic Varieties") et du r\'eseau europ\'een HPRN-CT-2002-00287.}

\bigskip

 {\bf \S 1. Rappels et pr\'eliminaires}

\bigskip

La cohomologie employ\'ee est la cohomologie \'etale, qui, sur un corps, donne la cohomologie
galoisienne. On supposera le lecteur familier avec la th\'eorie fine des tores alg\'ebriques
([CT/San1], [Vos2], [Vos3]), en particulier avec les notions de tore quasitrivial, de tore flasque, et
avec les divers types de r\'esolutions flasques (voir aussi [CT/San3]). On utilisera   librement
des r\'esultats sur les groupes lin\'eaires connexes quelconques que l'on peut trouver dans
[San]. On utilisera   librement la notion de $R$-\'equivalence sur les points $k$-rationnels
d'une vari\'et\'e alg\'ebrique d\'efinie sur un corps $k$ (voir [CT/San1], [Vos3]). Enfin
 on utilisera librement la cons\'equence suivante du th\'eor\`eme d'Hironaka  : si $k$ est un
corps de caract\'eristique z\'ero et $f : X\to Y$ est un morphisme de $k$-vari\'et\'es
quasiprojectives lisses int\`egres, il existe un morphisme $f_c : X_c \to Y_c$ de
$k$-vari\'et\'es projectives lisses int\`egres \'etendant $f$ (cf. [B/K1], 1.2.2).

\medskip

{\bf Proposition 1.1}  {\it Soient $k$ un corps,
$\k$ une cl\^oture s\'eparable de $k$ et $g=\gal(\k/k)$. Soit $X$ une $k$-vari\'et\'e
g\'eom\'etriquement int\`egre telle que $\k^{\times} = \k[X]^{\times}$, o\`u
$\k[X]^{\times}$ est le groupe des fonctions inversibles sur ${\overline X}$. La suite spectrale
de Leray pour le faisceau
$\G_m$ et la projection $X
\to {\rm Spec}(k)$ donne  naissance \`a une suite exacte naturelle
$$ 0 \to \pic(X) \to \pic({\overline X})^g \to
\br(k) \to \br_1(X) 
\to H^1(g,\pic({\overline X})) \to H^3(g,\k^{\times}),$$
o\`u  $\br_1(X)= {\rm Ker} [\br(X) \to \br({\overline X})]$.
 La fl\`eche
$H^1(g,\pic({\overline X})) \to H^3(g,\k^{\times})$ est nulle dans chacun des cas suivants :  $
X(k) \neq \emptyset$; $k$ est un corps $p$-adique ou r\'eel; $k$ est un corps de nombres. }

\medskip {\it D\'emonstration} C'est une cons\'equence bien connue de la suite spectrale de
Leray pour la topologie
\'etale, le morphisme $X \to {\rm Spec}(k)$ et le faisceau $\G_m$, et de la nullit\'e de
$H^3(g,\k^{\times})$ pour chacun des corps cit\'es.
\cqfd

\medskip

Soient $k$ un corps et  $G$ un $k$-groupe lin\'eaire. On note ${\hat G}={\rm Hom}_{{\overline
k}-{\rm groupes}}({\overline G}, {\Bbb G}_{m,{\overline k}})$ le groupe des caract\`eres
de $G$.  C'est un $g$-module continu discret de type fini; si $G$ est connexe,
ce module n'a pas de $\Z$-torsion.

\medskip

Soient $k$ un corps, $H$  un $k$-groupe alg\'ebrique et $Y$ une $k$-vari\'et\'e alg\'ebrique.
Un torseur (\`a gauche) sur $Y$ sous $H$ est une $k$-vari\'et\'e $X$ \'equip\'ee d'un morphisme
fid\`element plat $p : X  \to Y$ et d'une action $H \times_k X \to X$ du groupe $H$ sur $X$,
respectant le morphisme $p$, et telle que
le morphisme $H \times_k X \to X \times_YX$ donn\'e par $(h,x)  \mapsto (hx,x)$ soit un isomorphisme.
On a la d\'efinition analogue pour un torseur \`a droite. On dit parfois espace principal
homog\`ene au lieu de torseur.
\medskip

{\bf Proposition 1.2} (Sansuc [San], Prop. 6.10) {\it Soient $k$ un corps de caract\'eristique
z\'ero, $H$ un
$k$-groupe lin\'eaire connexe, $Y$ une $k$-vari\'et\'e lisse g\'eom\'etriquement connexe et $f :
X \to Y$ un torseur sur $Y$ sous $H$. On a alors une suite exacte
$$ 0 \to \k[Y]^{\times} \to \k[X]^{\times} \to {\hat H} \to
\pic({\overline Y}) \to \pic({\overline X}).  $$  }

Rappelons ici que le groupe des caract\`eres $ {\hat H}$ du groupe connexe $H$
s'identifie au module galoisien ${\overline k}[H]^{\times}/{\overline k}^{\times}$ (lemme de
Rosenlicht).
\medskip
{\bf Proposition 1.3}Ê {\it   Soient $k$ un corps de caract\'eristique z\'ero et
 $X/k$ une $k$-vari\'et\'e projective, lisse, g\'eom\'etriquement connexe. Soit $F/k$ une
extension de corps, avec $k$ alg\'ebriquement ferm\'e dans $F$. Soient ${\overline k} \subset
{\overline F}$ des cl\^otures alg\'ebriques de $k$ et $F$. Soient
$g_k={\rm Gal}({\overline k}/k)$ et $g_F={\rm Gal}({\overline F}/F)$. Supposons en outre que le
groupe ${\rm Pic}(X \times_k {\overline k})$ est un groupe ab\'elien de type fini sans torsion.
Alors :

{\rm (i) } Le groupe de Galois de ${\overline F}$ sur le corps compos\'e
$L={\overline k}.F$ agit trivialement sur le $g_F$-module ${\rm Pic}(X
\times_k {\overline F})$, qui est donc un $g$-module. La fl\`eche naturelle
${\rm Pic}(X \times_k {\overline k}) \to {\rm Pic}(X \times_k {\overline F})$ est un
isomorphisme de $g$-modules.

{\rm (ii)}  Chacune des propri\'et\'es suivantes vaut pour le $g_k$-module ${\rm Pic}(X \times_k
{\overline k})$ si et seulement si elle vaut pour le $g_F$-module
${\rm Pic}(X \times_k {\overline F})$ : \^etre un module de permutation, \^etre un facteur
direct d'un module de permutation,  
\^etre un module flasque, \^etre \'egal \`a un  $g$-r\'eseau donn\'e \`a addition pr\`es
de modules de permutation.

{\rm (iii)}  Les groupes $H^1(g_k,{\rm Pic}(X \times_k {\overline k}))$ et $H^1(g_F,{\rm Pic}(X \times_k
{\overline F}))$ sont isomorphes.

{\rm (iv) } Le quotient  de
$\br_1(X)$  par l'image de $\br(k)$ s'injecte dans le
quotient de $\br_1(X_F)$ par l'image de $\br(F)$. }

\medskip

{\it D\'emonstration} L'hypoth\`ese que ${\rm Pic}(X \times_k {\overline k})$  est un groupe de
type fini sans torsion \'equivaut \`a la combinaison de deux faits  : le groupe de cohomologie
coh\'erente
$H^1(X,{\cal O}_X)$ est nul et le groupe de N\'eron-Severi de $X
\times_k{\overline k}$ est sans torsion. La nullit\'e de $H^1(X,{\cal O}_X)$
\'equivaut au fait que la vari\'et\'e de Picard de $X$ est triviale. Le groupe de  Picard de $X
\times_k {\overline k}$ co\"{\i}ncide donc avec le groupe de N\'eron-Severi de $X \times_k
{\overline k}$, de m\^eme pour $X \times_k {\overline  F}$, et il est bien connu que le groupe de
N\'eron-Severi ne change pas par extension de corps de base alg\'ebriquement clos (ceci se
d\'eduit facilement du th\'eor\`eme analogue pour la cohomologie \'etale \`a coefficients de
torsion). Ainsi les injections naturelles
$$ {\rm Pic}(X \times_k {\overline k}) \hookrightarrow  {\rm Pic}(X \times_k  L)
\hookrightarrow  {\rm Pic}(X \times_k {\overline F})$$ sont des isomorphismes.
 Ceci \'etablit
l'\'enonc\'e (i), qui implique imm\'ediatement l'\'enonc\'e (ii) (si
$M$ est un ${\rm Gal}({\overline F}/F)$-module de permutation, alors
le groupe des points fixes de $M$ sous ${\rm Gal}({\overline F}/L)$ est
un ${\rm Gal}(L/F)$-module de permutation, c'est-\`a-dire un $g$-module de permutation).
L'\'enonc\'e (iii)
r\'esulte de la suite de restriction-inflation et de la nullit\'e de
$H^1(h,M)$ pour $M$ un groupe ab\'elien sans torsion
\'equip\'e d'une action triviale d'un groupe $h$ (profini). L'\'enonc\'e  (iv) r\'esulte de
(iii) et de la proposition 1.1. \cqfd

\bigskip

{\it Remarque 1.3.1} Par un r\'esultat de Serre, l'hypoth\`ese  que ${\rm Pic}(X \times_k
{\overline k})$ est un groupe de type fini sans torsion est satisfaite si la
${\overline k}$-vari\'et\'e $X
\times_k{\overline k}$ est unirationnelle, ce qui est certainement le cas pour une
compactification lisse d'un espace homog\`ene d'un groupe lin\'eaire connexe.

\bigskip

{\it Remarque 1.3.2}  L'int\'er\^et de la proposition 1.3 est que, pour
\'etablir  l'une des propri\'et\'es voulues pour $X/k$, on peut supposer
$X(k)\neq \emptyset$. Par la proposition ci-dessus, il suffit en effet 
d'\'etudier $X\times_kF/F$, o\`u $F=k(X)$ est le corps des fonctions de $X$.
La $F$-vari\'et\'e $X\times_kF$ poss\`ede un point $F$-rationnel, donn\'e par le
point g\'en\'erique de $X/k$. C'est ``l'astuce du
passage au point g\'en\'erique". Il existe d'autres variantes. Par exemple dans [B/K2],  on
consid\`ere une compactification lisse $X$ d'un groupe alg\'ebrique r\'eductif connexe $G$ non
n\'ecessairement quasid\'eploy\'e, et l'on passe au corps  des fonctions de la vari\'et\'e des
sous-groupes de Borel de $G$, ce qui a pour effet de quasid\'eployer $G$. 
\medskip

Pour le th\'eor\`eme suivant, on renvoie \`a [Bog] et [CT/San4] (Thm. 9.13).

\medskip {\bf Th\'eor\`eme 1.4} (Bogomolov) {\it Soient $k$ un corps alg\'ebriquement clos de
caract\'eristique z\'ero, $G$ un
$k$-groupe semi-simple simplement connexe et $H \subset G$ un sous-groupe ferm\'e connexe. Soit
$X_c$ une compactification lisse du quotient $X=G/H$. Alors $\br(X_c)=0$. \cqfd}

On peut se demander si ce th\'eor\`eme vaut encore si $G$ est une extension d'un tore
par un groupe semi-simple simplement connexe. Si tel \'etait le cas, on pourrait dans 
le th\'eor\`eme A remplacer $\br_1(X_c)$ par $\br(X_c)$.

\bigskip

Nous donnons maintenant une d\'efinition du {\it $k$-tore 
associ\'e} \`a un espace homog\`ene du type consid\'er\'e dans cet article.
Soient $G$ un $k$-groupe lin\'eaire connexe
et $X/k$ un espace homog\`ene de $G$
de stabilisateur g\'eom\'etrique ${\overline H}$ connexe.
Comme $X$ est
g\'eom\'etriquement int\`egre, le corps $k$ est alg\'ebriquement ferm\'e dans le corps des
fonctions $k(X)$ de $X$, donc $\k(X)=k(X)\otimes_k\k$ est le corps des fonctions de ${\overline
X}$, et le groupe de Galois de
$\k(X)$ sur $k(X)$ s'identifie au groupe de Galois $g$ de $\k$ sur $k$. Si l'on \'etend le corps
de base de $k$ \`a $k(X)$, le point g\'en\'erique de $X$ permet d'identifier $X_{k(X)}=X
\times_kk(X)$ au quotient de $G_{k(X)}$ par un $k(X)$-sous-groupe ferm\'e connexe $H_0$, groupe
qui sur
$\k(X)$ est isomorphe \`a ${\overline H}\times_{\k}\k(X)$. Si l'on prend le quotient du
$k(X)$-groupe  $H_0$ par son radical unipotent, puis le quotient du groupe obtenu par son groupe
d\'eriv\'e, on obtient un
$k(X)$-tore $T_0$ qui sur $\k(X)$ vient de $\k$, et qui est donc d\'eploy\'e par l'extension
$\k(X)$ sur $k(X)$. Le groupe des caract\`eres  de $T_0$ est donc de fa\c con naturelle un
$g$-r\'eseau, donc le dual d\'efinit un unique $k$-tore $T$ tel que $T\times_kk(X) \simeq T_0$.
On dira que $T$ est {\it le $k$-tore  associ\'e au $G$-espace homog\`ene $X$}.
 Le $g$-module ${\hat T}$ des caract\`eres de
$\overline T$ est  aussi le groupe des caract\`eres ${\hat H}$  de ${\overline H}$. Lorsque $X$
poss\`ede un point $k$-rationnel, le choix d'un
$k$-point $x \in X(k)$ d\'etermine un $k$-groupe alg\'ebrique
$H_x$, le fixateur de $x$, satisfaisant $H_x \times_k\k \simeq {\overline H}$. Si l'on quotiente
$H_x$ par son radical unipotent, puis le quotient obtenu par son groupe d\'eriv\'e, on obtient
un $k$-tore
$H_x^{\rm tor}$, le plus grand quotient torique de $H_x$, qui est isomorphe au $k$-tore $T$
associ\'e au $G$-espace homog\`ene $X=G/H_x$.

\bigskip

Terminons ce paragraphe en rappelant la d\'efinition et quelques propri\'et\'es des
groupes quasitriviaux  ([CT]). Les r\' esultats annonc\' es dans [CT] le sont pour les groupes
r\' eductifs connexes, le cas des groupes lin\' eaires connexes
quelconques s'en d\'eduit ais\'ement, en caract\'eristique nulle.
Soit $k$ un corps de caract\'eristique nulle, $\overline k$ une cl\^oture alg\'ebrique,
$g={\rm Gal}({\overline k}/k)$. Soit $G$ un $k$-groupe lin\'eaire connexe.
Le quotient de $G$ par son radical unipotent est un $k$-groupe r\'eductif connexe $G^{\rm red}$.
Le groupe $G^{\rm red}$ est une extension d'un $k$-tore $G^{\rm tor}$ par le groupe d\'eriv\'e
$G^{\rm ss}$ de $G^{\rm red}$, qui est un $k$-groupe semi-simple connexe.
Le $k$-groupe $G$ est dit quasitrivial si les deux propri\'et\'es suivantes sont satisfaites :

\smallskip

(i) Le $k$-tore $G^{\rm tor}$ est  quasitrivial (son groupe des caract\`eres est un $g$-module
de permutation).

(ii) Le groupe semi-simple $G^{\rm ss}$ est simplement connexe.
\smallskip

\noindent Ceci \'equivaut  ([CT], Prop. 1.11) \`a la combinaison des deux propri\'et\'es suivantes :
\smallskip

(iii) ${\hat{G}}  \simeq {\overline k}[G]^{\times}/{\overline k}^{\times} $ est un $g$-module de permutation.

(iv)  ${\rm Pic}({\overline G})=0.$

\smallskip

Pour tout $k$-groupe lin\'eaire connexe $G$, il existe une suite exacte 
$$ 1 \to F \to G_1 \to G \to 1$$
o\`u $G_1$ est un $k$-groupe lin\'eaire connexe quasitrivial et $F$ est
un $k$-tore flasque (de groupe des caract\`eres un $g$-module flasque)
central dans $G_1$. 

\medskip

{\bf Lemme 1.5} {\it Soit $k$ un corps de caract\'eristique nulle.
Soit $G$ un $k$-groupe lin\'eaire connexe et $X$ une $k$-vari\'et\'e qui est
un espace homog\`ene de $G$ \`a stabilisateur g\'eom\'etrique connexe.
Il existe un $k$-groupe alg\'ebrique lin\'eaire connexe quasitrivial $G_1$
tel que $X$ est un espace homog\`ene de $G_1$ \`a stabilisateur connexe.}
\medskip
{\it D\'emonstration} 
Soit ${\overline H}$ le stabilisateur g\'eom\'etrique pour l'action de $G$ sur $X$.
Soit 
$$1 \to F \to G_1 \to  G \to 1$$
une suite exacte comme ci-dessus. La $k$-vari\'et\'e $X$
est un espace homog\`ene de $G_1$. Le stabilisateur g\'eom\'etrique
pour l'action de $G_1$ sur $X$ est une extension  (centrale) de 
${\overline H}$ par le ${\overline k}$-tore ${\overline F}$.
C'est un groupe connexe. \cqfd

\bigskip

{\bf \S 2  Quotients de tores} 

\bigskip

{\bf Proposition 2.1} {\it Soit $k$ un corps de caract\'eristique z\'ero. Soit $ 1 \to T \to P
\to Q \to 1$ une suite exacte de $k$-tores, avec $P$ quasitrivial. Soit $Q_c$ une
$k$-compactification lisse du $k$-tore $Q$. Soit $ 1 \to T \to F \to P_1 \to 1$ une suite exacte
de $k$-tores avec $F$ flasque et $P_1$ quasitrivial.

{\rm (i) } Les $g$-modules $\Z$-libres de type fini ${\hat F}$ et
$\pic({\overline Q}_c)$ sont isomorphes \`a addition pr\`es de modules de permutation, en
particulier $\pic({\overline Q}_c)$ est un $g$-module flasque.

{\rm (ii)}  On a $$\br(Q_c)/\br(k)=H^1(g, \pic({\overline Q}_c)) =H^1(k,{\hat F})=\X^1_{\omega}(k,{\hat
T})=\X^2_{\omega}(k,{\hat Q}).$$

{\rm (iii) } $Q_c(k)/R= Q(k)/R=H^1(k,F)$.}

\medskip {\it D\'emonstration} D'apr\`es ([CT/San3], (1.3.2)), pour tout $k$-tore $T$ il existe
une suite exacte $ 1 \to T \to F \to P_1 \to 1$  du type indiqu\'e, et le $k$-tore flasque $F$
est bien d\'etermin\'e \`a multiplication par un $k$-tore quasitrivial pr\`es.
Formons l'expuls\'e (en
anglais : ``push-out") de
$T
\to P$ et $T \to F$. Ceci donne naissance \`a un diagramme commutatif de suites exactes de
$k$-tores :
$$
\diagram{ && 1 &&  1  &&  1  &&  &\cr &&\vfl{}{}{4mm}&&\vfl{}{}{4mm}&&  \vfl{}{}{4mm} &\cr 1
&\hfl{}{}{4mm}&T &\hfl{}{}{4mm}&P &\hfl{}{}{4mm}& Q &\hfl{}{}{4mm}& 1\cr
&&\vfl{}{}{4mm}&&\vfl{}{}{4mm}&&\pafl{}{}{4mm}&\cr 1 &\hfl{}{}{4mm}& F &\hfl{}{}{4mm}& N
&\hfl{}{}{4mm}&
 Q &\hfl{}{}{4mm}& 1\cr &&\vfl{}{}{4mm}&&\vfl{}{}{4mm}&&   &\cr
   &&P_1 & = & P_1 && && \cr &&\vfl{}{}{4mm}&&\vfl{}{}{4mm}&&   &\cr &&1 &&1&&   &\cr }
$$ Comme les tores $P$ et $P_1$ sont quasitriviaux, la suite m\'ediane verticale est scind\'ee,
d'o\`u une suite exacte
$$ 1\to F\to P\times P_1\to Q\to 1.
$$ Cette derni\`ere suite est une r\'esolution flasque du tore $Q$.  
Un r\'esultat de Voskresenski\u{\i} (cf. [CT/San1], Prop.~6) donne alors (i). 
La proposition 1.1 donne la premi\`ere \'egalit\'e dans
(ii). La seconde r\'esulte de (i) et de l'annulation de $H^1(g,{\hat P})$ pour un
$g$-module de permutation ${\hat P}$. Les deux derni\`eres \'egalit\'es  de (ii) s'\'etablissent
en consid\'erant les suites exactes de caract\`eres associ\'ees aux suites exactes de tores de
l'\'enonc\'e, et en utilisant   l'annulation mentionn\'ee \`a l'instant, l'annulation
de $H^1(h,{\hat F})$ pour ${\hat F}$
un $g$-module flasque et $h$ un sous-groupe ferm\'e procyclique de $g$, et
l'annulation de $\X^2_{\omega}(k,{\hat P})$ pour ${\hat P}$ un $g$-module de permutation. Quant
\`a (iii), c'est une application du Th\'eor\`eme 2  p. 199 et de la Proposition 13 p. 203 de
[CT/San1].
\cqfd

\bigskip

{\bf  \S 3. Espaces homog\`enes de groupes lin\'eaires connexes dont le stabilisateur
g\'eom\'etrique n'a pas de quotient torique}

\bigskip

{\bf Proposition 3.1} {\it Soient $k$ un corps de caract\'eristique z\'ero, $\k$ une cl\^oture
alg\'ebrique,
$g=\gal(\k/k)$,
$G$ un $k$-groupe quasitrivial, $X$ un $k$-espace homog\`ene sous $G$
\`a stabilisateur g\'eom\'etrique  $\overline H$ connexe. Supposons 
${\overline H}^{\rm tor}=1$.
Soit $X_c$ une $k$-compactification lisse de $X$. Alors :

{\rm (i)}  Le $g$-module $\pic({\overline X}_c)$ est stablement de permutation.

{\rm (ii) } La fl\`eche naturelle $\br(k) \to \br_1(X_c) $ est 
un isomorphisme.}

\medskip {\it D\'emonstration}  Pour \'etablir le point (i), on peut, d'apr\`es la
proposition 1.3 et la remarque 1.3.2, 
supposer $X(k) \neq \emptyset$, et donc $X=G/H$ avec $H$
un $k$-groupe lin\'eaire connexe tel que $H^{\rm tor}=1$. D'apr\`es le th\'eor\`eme de
Hironaka, il existe une $k$-compactification lisse $G_c$ de $G$ et un $k$-morphisme
$p : G_c \to X_c$ \'etendant le $k$-morphisme naturel $p : G \to G/H$.
 Comme $G$ est quasitrivial,  $\pic({\overline G})=0$, et 
${\overline k}[G]^{\times}/{\overline k}^{\times}=\hat{G}$
est un $g$-module de permutation. Sous l'hypoth\`ese ${\overline H}^{\rm tor}=1$,
on d\'eduit alors de la proposition 1.2 l'\'egalit\'e $\pic({\overline X})=0$
et le fait que la fl\`eche naturelle 
${\overline k}[X]^{\times}/{\overline k}^{\times} \to {\overline k}[G]^{\times}/{\overline k}^{\times}$
est un isomorphisme. Le $g$-module ${\overline k}[X]^{\times}/{\overline k}^{\times}$
est en particulier un $g$-module de permutation.
On dispose alors du diagramme commutatif de  suites exactes courtes 
 $$\diagram{0 & \hfl{}{}{4mm} & {\overline k}[G]^{\times}/{\overline k}^{\times} & \hfl{}{}{4mm} &
{\rm Div}_{\infty}({\overline G}_c) &\hfl{}{}{4mm} & \pic({\overline G}_c) & \hfl{}{}{4mm} & 0
\cr
&  & \ufl{}{f^*}{4mm} &  &
\ufl{}{f^*}{4mm} & & \ufl{}{f^*}{4mm} &  & 
\cr
0 & \hfl{}{}{4mm} & {\overline k}[X]^{\times}/{\overline k}^{\times} & \hfl{}{}{4mm} &
{\rm Div}_{\infty}({\overline X}_c) &\hfl{}{}{4mm} & \pic({\overline X}_c) & \hfl{}{}{4mm} & 0.
\cr
}
$$
Dans ce diagramme, la fl\`eche verticale de gauche est ici un isomorphisme.
C'est un th\'eor\`eme g\'en\'eral de [B/K2] que le $g$-module $\pic({\overline G}_c)$
est un $g$-module flasque. Comme ${\overline k}[G]^{\times}/{\overline k}^{\times}$ est ici un 
$g$-module de permutation, la suite horizontale sup\'erieure est donc scind\'ee.
La suite exacte inf\'erieure est la r\'etrotirette (en anglais : ``pull-back") de la suite
sup\'erieure par l'application $f^* : \pic({\overline X}_c) \to \pic({\overline G}_c)$. Elle est
donc aussi scind\'ee. Comme ${\overline k}[X]^{\times}/{\overline k}^{\times}$ et 
${\rm Div}_{\infty}({\overline X}_c)$ sont tous deux des $g$-modules de permutation,
ceci \'etablit le point (i).

Etablissons le point (ii). Lorsque $X(k) \neq \emptyset$, on a vu ci-dessus
que les modules galoisiens ${\overline k}[X]^{\times}/{\overline k}^{\times}$
et $ {\overline k}[G]^{\times}/{\overline k}^{\times}$ sont isomorphes.
Montrons comment l'astuce du passage au point g\'en\'erique permet de montrer qu'il en
est encore ainsi sans supposer $X(k) \neq \emptyset$. 
Soit $E=k(X)$ le corps des fonctions de $X$
et $L={\overline k}(X)$ le corps des fonctions de ${\overline X}$.
L'extension $L/E$ est galoisienne de groupe $g$.
C'est un fait g\'en\'eral que pour tout corps alg\'ebriquement clos $k_1$,
toute extension de corps
 $k_1 \subset k_2$ et toute $k_1$-vari\'et\'e lisse
connexe $X$, la fl\`eche naturelle $k_1[X]^{\times}/k_1^{\times} \to k_2[X]^{\times}/k_2^{\times}$
est un isomorphisme. Les fl\`eches ${\overline k}[X]^{\times}/{\overline k}^{\times} \to
L[X]^{\times}/L^{\times}$ et ${\overline k}[G]^{\times}/{\overline k}^{\times} \to
L[G]^{\times}/L^{\times}$, qui sont des homomorphismes de $g$-modules, sont donc
des isomorphismes. La $E$-vari\'et\'e $X\times_KE$ poss\`ede un 
point rationnel, correspondant au point g\'en\'erique de $X$.
L'argument donn\'e au point (i) montre alors que
l'on dispose d'un $g$-isomorphisme $ L[X]^{\times}/L^{\times} \simeq L[G]^{\times}/L^{\times}$.
Ainsi les $g$-modules ${\overline k}[X]^{\times}/{\overline k}^{\times}$
et $ {\overline k}[G]^{\times}/{\overline k}^{\times}$ sont isomorphes.
En particulier ${\overline k}[X]^{\times}/{\overline k}^{\times}$
est un $g$-module de permutation. Le th\'eor\`eme 90 de Hilbert
assure alors que la suite exacte \'evidente
$$ 1 \to {\overline k}^{\times} \to {\overline k}[X]^{\times} \to {\overline k}[X]^{\times}/{\overline k}^{\times}
\to 1$$
est scind\'ee. L'application naturelle $H^2(g,{\overline k}^{\times}) \to H^2(g,{\overline k}[X]^{\times})$
est donc injective. Comme on a $\pic({\overline X})=0$, la suite spectrale de Leray pour
le morphisme $X \to {\rm Spec}(k)$ et le faisceau \'etale $\G_m$ donne
$H^2(g,{\overline k}[X]^{\times}) \simeq {\rm Ker}[\br(X) \to \br({\overline X})].$
Ainsi l'application naturelle $\br(k) \to \br(X)$ est injective,  a fortiori en est-il
de m\^eme de l'application
$\br(k) \to \br(X_c)$. D'apr\`es le point (i), le $g$-module $\pic({\overline X}_c)$
est stablement de permutation, en particulier $H^1(g,\pic({\overline X}_c))=0$.
Une application de la proposition 1.1 \`a $X_c$ donne alors le point (ii).

\medskip

\bigskip

{\bf \S 4. Le th\'eor\`eme B}

\bigskip

L'\'enonc\'e suivant rassemble des propri\'et\'es connues.

\medskip

{\bf Th\'eor\`eme 4.1} {\it Soit $k$ un corps de caract\'eristique z\'ero, de
dimension cohomologique 1. Le corps $K=k((t))$ des s\'eries formelles en une variable sur
$k$ satisfait les
propri\'et\'es suivantes.

{\rm (i)}  Sa dimension cohomologique est 2.

{\rm (ii) } Sur toute extension finie de $K$, indice et exposant
des alg\`ebres simples centrales co\"{\i}ncident.

{\rm (iii)} Pour tout $K$-groupe quasitrivial $G$
sur $K$, on a $H^1(K,G)=1$.

{\rm (iv)} Soient $G$ un $K$-groupe semi-simple simplement connexe,
$\mu$ son centre et $$1 \to \mu \to G \to G^{\rm ad} \to 1$$
l'isog\'enie associ\'ee. Alors la fl\`eche de bord
$H^1(K,G^{\rm ad}) \to H^2(K,\mu)$ est une bijection.

{\rm (v) } Soient ${\overline H}$ un ${\overline K}$-groupe lin\'eaire connexe
et $L=({\overline H},\kappa)$ un $K$-lien. Soit $H^{\rm tor}$
le $K$-tore associ\'e \`a $L$. On dispose d'une application naturelle
d'ensembles $H^2(K,L) \to H^2(K,H^{\rm tor})$.
 Un \'el\'ement $\eta \in H^2(K,L)$
est neutre si et seulement si son image dans $H^2(K,H^{\rm tor})$
est triviale.

(vi) Soient $G$ un $K$-groupe quasitrivial
et $E$ un espace homog\`ene de $G$. Si le 
stabilisateur g\'eom\'etrique
est connexe et n'a pas de quotient torique, alors $E(K)\neq \emptyset$.}

\medskip

{\it Indications}  Les \'enonc\'es (i), (ii) et (iii) sont bien connus. On
trouvera des r\'ef\'erences pr\'ecises dans [CTGiPa] (Thm. 1.5).
 L'\'enonc\'e (iii) est d\^u \`a Bruhat et Tits lorsque $G$
est un $K$-groupe semi-simple simplement connexe, et le
cas g\'en\'eral r\'esulte du th\'eor\`eme 90 de Hilbert
et de la cohomologie galoisienne de la suite exacte
$1 \to G^{\rm ss} \to G \to G^{\rm tor} \to 1. $
Les \'enonc\'es (iv) et (v) sont
dans la th\`ese de Douai, lequel s'appuie sur les r\'esultats  de Bruhat
et Tits. Comme
\'etabli  dans [CTGiPa] (Thm. 2.1 (a)), l'\'enonc\'e (iv) vaut pour tout
corps $K$ de caract\'eristique nulle satisfaisant les
\'enonc\'es (i), (ii) et (iii). De m\^eme, l'\'enonc\'e (v)
vaut pour tout tel corps : c'est le Th\'eor\`eme 5.4 de [CTGiPa],
dont la d\'emonstration combine la m\'ethode de Borovoi ([Bo1], cas
des corps $p$-adiques et des corps globaux) et l'\'enonc\'e (iv).
Comme rappel\'e par Borovoi [Bo1],
 la donn\'ee de $(G,E)$
dans (vi) d\'efinit un $K$-lien $L=({\overline H},\kappa)$ et le $K$-tore
associ\'e $H^{\rm tor}$ est par hypoth\`ese trivial. L'\'enonc\'e (vi)
r\'esulte alors de (v). \cqfd

\bigskip

{\it Remarque 4.1.1} Les indications donn\'ees montrent que le th\'eor\`eme 4.1
vaut plus g\'en\'eralement sur un corps $K$ de caract\'eristique nulle
satisfaisant les conditions (i), (ii) et (iii).

On appellera ici {\it bon corps de dimension cohomologique $\leq 2$}
un corps de caract\'eristique nulle satisfaisant ces trois conditions.
Des exemples de tels corps sont les corps de nombres totalement imaginaires,
les corps $p$-adiques, les corps de fractions d'anneaux locaux hens\'eliens
de dimension 2 \`a corps r\'esiduel alg\'ebriquement clos, 
les corps de fonctions de deux variables sur les
complexes. Pour ces derniers corps l'\'enonc\'e (ii) est un  r\'esultat
 de de Jong et l'\'enonc\'e (iii) vaut en l'absence
de type $E_8$. On trouvera des r\'ef\'erences dans [CTGiPa].
\bigskip

Comme indiqu\'e dans l'introduction, nous devons \`a O. Gabber 
l'id\'ee de la d\'emonstration du th\'eor\`eme suivant.

\medskip

{\bf Th\'eor\`eme 4.2} (Th\'eor\`eme B) {\it Soit $A$ un anneau de valuation
discr\`ete de corps des fractions $K$, de corps r\'esiduel $k$ de
caract\'eristique nulle. Soient $G$ un $K$-groupe quasitrivial et 
$E/K$ un espace homog\`ene sous $G$ de stabilisateur g\'eom\'etrique un
groupe connexe sans quotient torique. Soit $X$ un
$A$-sch\'ema propre, r\'egulier, int\`egre, dont la fibre g\'en\'erique
contient $E$ comme ouvert dense. Alors il existe une composante de
multiplicit\'e 1 de la fibre sp\'eciale de $X/A$ qui est g\'eom\'etriquement
int\`egre sur son corps de base $k$.}

\medskip

{\it D\'emonstration}  Un proc\'ed\'e classique, qui appara\^{\i}t dans la d\'emonstration
du th\'eor\`eme de Merkur'ev et Suslin,  permet de
construire une extension $l/k$ avec $k$ alg\'ebriquement ferm\'e dans $l$
et $l$ de dimension cohomologique 1.  Indiquons le principe. S'il existe
 une extension finie $k_1$ de $k$ et
une vari\'et\'e de Severi-Brauer $W$ non triviale sur $k_1$, on plonge $k$
dans le corps des fonctions de la descendue \`a la Weil $R_{k_1/k}(W)$. On it\`ere
\`a l'infini le proc\'ed\'e. Pour plus de d\'etails, voir la note de Ducros [D].

Soit $\pi$ une uniformisante de $A$. Comme le corps r\'esiduel
 $k$ est de caract\'eristique z\'ero, le
compl\'et\'e de $A$ est isomorphe \`a $k[[t]]$ ([S], Chap. II, \S 4, Th\'eor\`eme 2),
 et l'on peut supposer que l'inclusion
$A \subset k[[t]]$ envoie $\pi $ sur $t$.
 D'apr\`es le th\'eor\`eme 4.1,
$E(l((t))) \neq \emptyset$. Donc $X(l((t))) \neq \emptyset$ et comme
$X/A$ est propre, $X(l[[t]]) \neq \emptyset$, c'est-\`a-dire qu'il
existe un morphisme $f : Spec(l[[t]]) \to X$ tel que le compos\'e
$Spec(l[[t]]) \to X \to Spec(A)$ est la fl\`eche d\'eduite
de $A \subset k[[t]] \subset B=l[[t]]$.

On consid\`ere le point (sch\'ematique)  $y \in X$ image du point ferm\'e de $Spec(l[[t]])$.
C'est un point de la fibre sp\'eciale de $X/A$. 
Soit $C$
l'anneau local de $X$ en $y$. Comme $X$ est r\'egulier, $C$ est un anneau local
r\'egulier, en particulier factoriel (Auslander-Buchsbaum, Serre).
 On a alors des
homomorphismes d'anneaux locaux $A \to C \to B$, l'homomorphisme compos\'e $A \to B$ n'\'etant
autre que l'inclusion \'evidente $A \subset k[[t]] \subset l[[t]]$. 

Soit $\pi_C =
u \cdot \prod_{i=1}^r
\rho_i^{n_i}$, avec les $n_i>0$, une d\'ecomposition de l'image de $\pi $ dans
$C$ : ici $u$ est une unit\'e de $C$, chaque $\rho_i$ appartient 
\`a l'id\'eal maximal de $C$ et est irr\'eductible. Cette d\'ecomposition correspond \`a la
description de la fibre sp\'eciale de $X$ au voisinage du point $y$. On a alors
l'\'egalit\'e 
$$t = u_B \cdot \prod_{i=1}^r (\rho_{i,B})^{n_i}\in B=l[[t]],$$ et l'image $\rho_{i,B}$
de  $\rho_i\in C$ dans $B=l[[t]]$ appartient
\`a l'id\'eal maximal de $B=l[[t]]$. 
 Ainsi $r=1$ et $n_1=1$, c'est-\`a-dire que $\pi_C$
est un \'el\'ement irr\'eductible de l'anneau local r\'egulier $C$. Ainsi la fibre sp\'eciale de
$X/A$ au voisinage du point $y$ n'a qu'une composante, de multiplicit\'e 1, soit $Z$.
Son corps r\'esiduel $k(Z)$ est le corps des fractions de $C/\pi_C$. L'homomorphisme compos\'e
$A \to C \to l[[t]]$ induit un homomorphisme compos\'e $A/\pi   \to C/\pi_C \to B/t$, soit encore
$k \to C/\pi_C \to l$. Comme $C/\pi_C$ est r\'egulier et en particulier normal,
la fermeture alg\'ebrique $k_1$ de $k$ dans le corps $k(Z)$ est un corps contenu dans $C/\pi_C$.
On a donc les inclusions de corps $k \subset k_1 \subset l$.
 Comme $k$ est alg\'ebriquement ferm\'e dans $l$, on  a $k=k_1$ et $k$ est alg\'ebriquement ferm\'e
dans $k(Z)$.
\cqfd
\bigskip

{\bf \S 5. Le th\'eor\`eme A}

\bigskip

{\bf Th\'eor\`eme 5.1} (Th\'eor\`eme A)  {\it Soient $k$ un corps de caract\'eristique nulle,
$G$ un $k$-groupe lin\'eaire connexe,
$X$ une $k$-vari\'et\'e espace homog\`ene de
$G$, de stabilisateur g\'eom\'etrique connexe.  Soit $X_c$ une $k$-compactification lisse de $X$.

{\rm (i) }  Le $g$-module
${\rm Pic}({\overline X}_c)$ est un $g$-module flasque, c'est-\`a-dire que pour tout sous-groupe
ferm\'e $h \subset g$, on a $H^1(h,\hom_\Z({\rm Pic}({\overline X}_c),\Z))=0$, soit encore
${\rm Ext}^1_h({\rm Pic}({\overline X}_c),\Z)=0$.

{\rm (ii)}  Pour tout sous-groupe ferm\'e procyclique $h \subset g$, on a  $H^1(h,{\rm Pic}({\overline
X}_c))=0$.

Supposons de plus $G$ quasitrivial, i.e. extension d'un $k$-tore
quasitrivial par un $k$-groupe simplement connexe. 

{\rm (iii) }   Soit 
$T$ le $k$-tore associ\'e au $G$-espace homog\`ene $X$. 
Soit $1 \to T \to F \to P \to 1$
une suite exacte de $k$-tores, avec $F$ flasque et $P$ quasitrivial. Les modules galoisiens
$\pic({\overline X}_c) $ et $\hat F$ sont \'egaux \`a addition pr\`es de modules de permutation.

{\rm (iv)}  Le quotient de $\br_1(X_c)$ 
 par l'image du groupe $\br(k)$ s'injecte 
dans le groupe
$\X^{1}_{\omega}(k,{\hat T})$, et est isomorphe \`a ce dernier groupe si $X(k)\neq \emptyset$ ou
si $k$ est un corps de nombres.}

\bigskip

{\it D\'emonstration} D'apr\`es le lemme 1.5, pour \'etablir ce th\'eor\`eme, on peut supposer d'embl\'ee le groupe $G$ quasitrivial.

D'apr\`es la proposition 1.3 et les remarques subs\'equentes, le
changement de base de $k$ au corps des fonctions de
$X$ permet de se ramener au cas o\`u $X(k) \neq \emptyset$, i.e. $X=G/H$ avec $H$ un
$k$-sous-groupe ferm\'e connexe de $G$.
Soit  $T=H^{\rm tor} $ et notons
$H_1$ le noyau de $H \to T$. C'est un $k$-groupe connexe extension d'un $k$-groupe semi-simple
par un $k$-groupe unipotent.  
Rappelons une construction g\'eom\'etrique   utilis\'ee par Borovoi (\S 4.2 de [Bo2]).
Soit
$$ 1 \to T \to P \to Q \to 1$$ une suite exacte de $k$-tores, avec $P$ quasitrivial.
Soit $Q_c$ une $k$-compactification  lisse du $k$-tore $Q$.
Le groupe $H$ agit sur $P$  via l'homomorphisme $H \to T$.
On d\'efinit une action  \`a droite de $H$ sur 
$ G \times P$ par la formule $(g,p) \cdot h=(gh,h^{-1}p)$.
Soit $Z=(G \times P)/H$ le quotient de $G \times P$ par cette action.

D'un c\^ot\'e $Z$ est un torseur sur $X$ sous le $k$-tore quasitrivial $P$, via l'action
\`a droite de $P$ sur le second facteur de $G \times P$. 
D'un autre c\^ot\'e $Z$
est un $G$-espace 
homog\`ene \`a gauche sur $Q$, via l'action \`a gauche de $G$ sur le
premier facteur de $G \times P$; un calcul simple montre que les
stabilisateurs g\'eom\'etriques pour cette action sont isomorphes au groupe
connexe ${\overline H}_1$, qui v\'erifie ${\overline H}_1^{\rm tor}=1$.

Le th\'eor\`eme d'Hironaka assure l'existence d'une  $k$-compactification  lisse
$Z_c$ de $Z$, d'un $k$-morphisme projectif   $f : Z_c \to Q_c$  \'etendant $Z \to Q=P/T$
et d'un $k$-morphisme projectif $Z_c \to X_c$ \'etendant $Z \to X$. 
 Comme $Z$ est un torseur sur
$X=G/H$ sous le $k$-tore quasitrivial $P$,  
le th\'eor\`eme 90 sous la forme de Grothendieck 
implique que
$Z$ est  $k$-birationnel au produit 
$X \times_kP$ de $X$  et de la $k$-vari\'et\'e $k$-rationnelle $P$.
Les modules galoisiens $\pic({\overline X}_c)$ et 
$\pic({\overline Z}_c)$ sont donc \'egaux \`a addition pr\`es de modules de permutation,
   $\br(X_c) =\br(Z_c)$ et $\br_{1}(X_c) =\br_{1}(Z_c)$. Il suffit
donc d'\'etablir le th\'eor\`eme en y rempla\c cant $X_c$ par $Z_c$.
Il existe un ouvert $U$ de $Q_c$, contenant $Q$, de compl\'ementaire dans $Q_c$ un ferm\'e $R$  de
codimension au moins 2, tel que toutes les fibres de $ f : V=f^{-1}(U) \to U$ soient
\'equidimensionnelles de dimension $\dim(Z) - \dim(Q)$. Soit $F_0=f^{-1}(R)  \subset Z_c$ le
ferm\'e compl\'ementaire de $V$.

Il y a un nombre fini de points $m$  de codimension 1 de $Q_c$, 
 dont la fibre $f^{-1}(m)$  n'est
pas g\'eom\'etriquement int\`egre sur le corps r\'esiduel $k(m)$. 
{\it D'apr\`es le th\'eor\`eme 4.2, on peut en chaque tel point  $m$ fixer une
composante de multiplicit\'e~1 g\'eom\'etriquement int\`egre de la fibre $f^{-1}(m)$}. Soit
$F_m \subset Z_c$ le ferm\'e  qui est la r\'eunion  des adh\'erences {\it des autres 
composantes} de $f^{-1}(m)$, et soit $F_{1} = \cup_{m} F_{m }$. C'est un ferm\'e propre de $Z_c$. Les points g\'en\'eriques des composantes de $F_{1}$,
qui sont de codimension 1 sur $Z_{c}$, n'appartiennent pas \`a $F_{0}$.
Soit $F=F_0 \cup F_1 \subset Z_{c}$.

On a   le diagramme de morphismes suivant :
$$
\diagram{ Z        &   \subset & V           &     \subset & Z_c \cr
      \vfl{}{}{4mm}&            &\vfl{}{}{4mm} &              &\vfl{}{f}{4mm} \cr 
        Q           &  \subset & U           & \subset       & Q_c.  \cr}
$$

Soient $K=k(Q)$ le corps des fonctions de $Q$ et $L={\overline k}(Q)={\overline k}.K$ 
le corps des
fonctions de ${\overline Q}$.
Soient $\eta={\rm Spec} (K)$  et ${\overline \eta}= {\rm Spec} (L)$.
Notons ${Z}_{c,{\overline \eta}}$
la fibre $Z_c \times_{Q_c} {\rm Spec}(L)$  de
${\overline f}$ au-dessus du point g\'en\'erique de ${\overline Q_c}$.
Soit  $M$ le $g$-module de permutation sur les points de codimension 1 de ${\overline Z}_c$ appartenant \`a ${\overline F}$.
La restriction de ${\overline Z}_c$ \`a ${Z}_{c,{\overline \eta}}$ donne naissance
\`a un complexe de $g$-modules :
 $$0 \to M  \oplus \pic({ \overline Q}_c) \to \pic({\overline Z}_c) \to 
 \pic( Z_{c,{\overline \eta}})
\to 0,   \hskip1cm  (5.1)$$
o\`u la fl\`eche $\pic({ \overline Q}_c) \to \pic({\overline Z}_c)$
est donn\'ee par $f^*$. 

\medskip

{\bf Lemme 5.1.1} {\it Ce complexe est une suite exacte}.
\medskip

{\it D\'emonstration}
La surjectivit\'e de $\pic({\overline Z}_c) \to 
 \pic( Z_{c,{\overline \eta}})$ r\'esulte de la lissit\'e de ${\overline Z}_c$.
 Montrons l'exactitude 
au terme m\'edian.  Soit $\Delta$ un diviseur de ${\overline Z}_c$
d'image nulle dans $ \pic( Z_{c,{\overline \eta}})$.
Le diviseur $\Delta$ est rationnellement \'equivalent sur ${\overline Z}_c$
\`a un diviseur dont la restriction \`a $Z_{c,{\overline \eta}}$ est nulle.
Pour \'etablir l'exactitude, on peut donc supposer que la restriction
de $\Delta$ \`a $ Z_{c,{\overline \eta}}$ est nulle, c'est-\`a-dire que $\Delta$
est une combinaison lin\'eaire \`a coefficients entiers de diviseurs 
irr\'eductibles (r\'eduits) $\delta $ sur ${\overline Z}_c$ tels que $f(\delta )$
soit de codimension au moins 1 dans ${\overline Q}_{c}$.
Si l'image $f(\delta )$ est de codimension au moins 2, alors $\delta $ est l'une des
composantes de $F_{0}$. Supposons $f(\delta )=\gamma$  de codimension 1.
Si l'image r\'eciproque de $\gamma$ par $f$ est un diviseur irr\'eductible de
${\overline Z}_c$ de multiplicit\'e 1, alors $\delta =f^*(\gamma)$ et
la classe de $\delta $ appartient \`a $f^*(\pic({ \overline Q}_c))$.
Si l'image r\'eciproque de $\gamma$ par $f$ n'est pas int\`egre, soit $\delta$ est une composante
de ${\overline F}_{1}$, soit la multiplicit\'e de $\delta $ dans $f^*(\gamma)$ est 1, et l'on a
$f^*(\gamma)= \delta  + \delta_{1}$, avec $\delta_{1}$ somme de composantes
de ${\overline F}_{1}$.  Dans tous les cas, la classe de $\delta $ dans $\pic({\overline Z}_c)$
appartient \`a l'image de $M  \oplus \pic({ \overline Q}_c)$.
Montrons l'exactitude
\`a gauche. Soit $\Delta_0$ un diviseur de ${\overline Z}_c$ \`a support
dans ${\overline F}_0$. Soit $\Delta_1$ un diviseur de ${\overline Z}_c$ \`a support
dans ${\overline F}_1$. Soit $\gamma$ un diviseur (non n\'ecessairement irr\'eductible) 
de ${\overline Q}_c$.
Supposons qu'il existe une fonction rationnelle $h \in {\overline k}(Z_c)^{\times}$
telle que
$$\Delta_0 + \Delta_1 + f^*(\gamma) =
{\rm div}_{{\overline Z}_c}(h) \in {\rm Div}({\overline Z}_c).$$
Comme la fibre g\'en\'erique de $f$ est projective, lisse,
et g\'eom\'etriquement int\`egre, cette \'egalit\'e implique que la fonction $h$ est
l'image r\'eciproque par $f^*$ d'une fonction que nous noterons encore $h \in
{\overline k}(Q_c)^{\times}$. On a alors $\Delta_0 + \Delta_1= f^*({\rm div}_{{\overline
Q}_c}(h) - \gamma) \in {\rm Div}({\overline Z}_c)$. La restriction de cette
\'egalit\'e \`a ${\overline V}$ se lit 
$\Delta_1= f^*({\rm div}_{{\overline U}}(h) - \gamma) \in {\rm Div}({\overline V})$. 
D'apr\`es la d\'efinition de $F_{1}$,  cette \'egalit\'e implique $\Delta_1=0 \in {\rm Div}({\overline Z}_{c})$ et 
${\rm div}_{\overline U}(h) - \gamma=0 \in {\rm Div}({\overline U})$. Puisque
le compl\'ementaire de $U$ dans la vari\'et\'e lisse $Q_c$ est de codimension au
moins 2, la derni\`ere \'egalit\'e implique
${\rm div}_{{\overline
Q}_c}(h) - \gamma=0 \in {\rm Div}({\overline Q}_c).$ 
On en d\'eduit $\Delta_0=0 \in {\rm Div}({\overline Z}_c)$, ce qui ach\`eve
d'\'etablir l'exactitude de la suite (5.1). \cqfd

Soit ${\overline K}$ une cl\^oture alg\'ebrique de $K$.
D'apr\`es la proposition 3.1,  la fl\`eche naturelle
$\br(L) \to \br({Z}_{c,{\overline \eta}}) $ est injective.
La proposition 1.1 assure alors que 
 le ${\rm Gal}({\overline K}/K)$-module
 $\pic( Z_{c,{\overline \eta}})$
est le module obtenu en prenant les points fixes de
$\pic(Z_c \times_{Q_c} {\rm Spec}({\overline K}))$
sous ${\rm Gal}({\overline K}/L)$.
D'apr\`es la proposition 3.1, le ${\rm Gal}({\overline K}/K)$-module
 $\pic(Z_c \times_{Q_c} {\rm Spec}({\overline K}))$ est
stablement de permutation, il en est donc de m\^eme du $g$-module
$\pic( Z_{c,{\overline \eta}})$.

D'apr\`es la proposition 3.1, la fl\`eche $\br(K) \to \br(Z_{c,{\eta}})$ est 
injective. Il en est donc de m\^eme 
de la fl\`eche  $f^* : \br(Q_c) \to \br(Z_c)$ (ce qui implique en particulier 
que la fl\`eche naturelle $\br(k) \to \br(Z_c)$ est une injection) puis de la fl\`eche
$\br(Q_c)/\br(k) \to \br(Z_c)/\br(k)$, et donc (d'apr\`es la proposition 1.1)
de $H^1(g,\pic({ \overline Q}_c)) \to H^1(g,\pic({\overline Z}_c))$.
Le m\^eme argument vaut en rempla\c cant $g$ par tout sous-groupe ferm\'e $h \subset g$.
Ainsi pour tout tel sous-groupe $h$, la suite exacte (5.1) induit 
une injection $H^1(h,\pic({ \overline Q}_c)) \to H^1(h,\pic({\overline Z}_c))$.

\medskip 
{\bf Lemme 5.1.2} {\it Soient $g$ un groupe profini et $0 \to A  \to B
\to C \to 0$ une suite exacte de $g$-modules continus de type fini comme groupes ab\'eliens.  
Supposons que $C$ est un $g$-module stablement de permutation et que pour
tout sous-groupe ferm\'e $h \subset g$, l'application induite $H^1(h,A)\to H^1(h,B)$ est 
injective (et donc un isomorphisme). Alors la suite est scind\'ee. } 
\medskip

{\it D\'emonstration} 
La suite exacte donn\'ee induit une suite exacte longue
$${\rm Hom}_g(C,B) \to {\rm Hom}_g(C,C) \to {\rm Ext}^1_g(C,A) \to {\rm Ext}^1_g(C,B).$$
Pour $h$ sous-groupe ferm\'e d'indice fini de $g$ et $M$ un $g$-module continu,
 on a un  isomorphisme
naturel (Shapiro)
${\rm Ext}^1_g(\Z[g/h],M)
\simeq H^1(h,M)$. L'hypoth\`ese assure donc que la fl\`eche ${\rm Ext}^1_g(C,A) \to {\rm
Ext}^1_g(C,B)$ est injective, l'application identit\'e dans ${\rm Hom}_g(C,C)$ se rel\`eve 
en un \'el\'ement de ${\rm Hom}_g(C,B)$, la suite est scind\'ee.
\cqfd

On voit alors que la suite exacte de $g$-modules (5.1) est scind\'ee,
ce qui montre que les  $g$-modules $\pic({\overline Z}_c)$ et
$\pic({ \overline Q}_c)$ sont \'egaux \`a addition pr\`es de $g$-modules
de permutation.  Il en est donc aussi de m\^eme des $g$-modules
$\pic({\overline X}_c)$ et
$\pic({ \overline Q}_c)$.

D'apr\`es Voskresenski\u{\i} (voir la proposition 2.1),
$\pic({ \overline Q}_c)$ est un $g$-module flasque. 
Le $g$-module $\pic({\overline X}_c)$ est donc flasque.
De plus sa classe \`a addition pr\`es d'un module de permutation
est donn\'ee par la proposition 2.1. Ceci \'etablit les points
(i) \`a (iii). Le point (iv) r\'esulte alors de la proposition 1.1.
\cqfd

\bigskip

{\bf \S 6. $R$-\'equivalence}

\bigskip

Dans ce paragraphe, nous utiliserons une propri\'et\'e fonctorielle simple des torseurs
universels ([CT/San2]). Soit $f  : X \to Y$ un morphisme de $k$-vari\'et\'es lisses,
projectives, g\'eom\'etriquement int\`egres. Soient $M \in X(k)$ et $N=f(M) \in Y(k)$. Supposons
que les modules galoisiens $ \pic({\overline X})$ et $ \pic({\overline Y})$ sont libres, de type
fini sur $\Z$. Soient $S_X$ et $S_Y$ les $k$-tores duaux de ces $g$-r\'eseaux. \`A la fl\`eche
$f^* : \pic({\overline Y}) \to 
\pic({\overline X})$ correspond un $k$-homomorphisme de $k$-tores
$S_X \to S_Y$. Soit ${\cal T}_X$ le torseur universel sur 
$X$ de fibre triviale en $M$ et soit ${\cal T}_Y$  le torseur universel sur 
$Y$ de fibre triviale en $N$. On a alors un isomorphisme de $S_Y$-torseurs sur $X$ :
$$ {\cal T}_X \times^{S_X} S_Y \simeq f^*({\cal T}_Y).$$ Ceci r\'esulte imm\'ediatement de la
d\'efinition des torseurs universels ([CT/San2], p. 408) et de la fonctorialit\'e de la suite
(2.0.2) (ibid.). Ceci implique que le diagramme suivant est commutatif
$$ \diagram{ X(k)  && \to && H^1(k,S_X) & \cr
                    \vfl{}{}{4mm} &&  && \vfl{}{}{4mm}  & \cr
          Y(k)  && \to && H^1(k,S_Y), &\cr}$$ o\`u les fl\`eches horizontales sont donn\'ees par
\'evaluation des torseurs universels ${\cal T}_X$ et ${\cal T}_Y$ sur les $k$-points. Rappelons
par ailleurs que les fl\`eches horizontales passent au quotient par la $R$-\'equivalence
([CT/San2], Prop. 2.7.2 p. 444).

\bigskip

 Soient $k$ un corps de caract\'eristique nulle, $X$ une $k$-vari\'et\'e
espace homog\`ene d'un $k$-groupe lin\'eaire connexe, de stabilisateur
g\'eom\'etrique connexe, et $e$ un $k$-point de $X$. Soient $X_c$ une
$k$-compactification lisse  de
$X$ et $F$ le $k$-tore de groupe des caract\`eres ${\hat F}={\rm
Pic}({\overline X}_c)$. Soit ${\cal T} \to  X_c$  le torseur universel sur $X_c$ de
fibre triviale en $e$.  C'est un torseur sur $X_c$ sous le $k$-tore $F$. 
On a une fl\`eche d'\'evaluation
associ\'ee $X_c(k) \to H^1(k,F)$, qui induit une fl\`eche
$ X_c(k)/R \to H^1(k,F).$ 

\medskip
{\bf Th\'eor\`eme 6.1} {\it 
Avec les notations et hypoth\`eses de l'alin\'ea
pr\'ec\'edent, si $k$ est
un bon corps  de dimension cohomologique
$\leq 2$ (cf. 4.1.1), alors la fl\`eche $ X_c(k)/R \to H^1(k,F)$  est surjective.}

\medskip

{\it D\'emonstration}  
D'apr\`es le lemme 1.5, il existe un $k$-groupe
quasitrivial $G$ sous lequel $X$ est un espace homog\`ene, tel que le 
stabilisateur de $e$ est un $k$-groupe lin\'eaire connexe.
 La construction g\'eom\'etrique utilis\'ee au   \S 5, dont on garde les notations,
permet de ramener la proposition pour $X$ \`a la proposition pour $Z$. On dispose 
alors de $f : Z_c \to Q_c$  \'etendant le $G$-espace homog\`ene $f : Z
\to Q$ dont tout stabilisateur  g\'eom\'etrique est isomorphe \`a un groupe connexe 
$\overline H$ satisfaisant ${\overline H}^{\rm tor}=1$.

On a le diagramme commutatif
$$ \diagram{ Z_c(k)/R && \to && H^1(k,F_{Z_c}) & \cr
                    \vfl{}{}{4mm} &&  && \vfl{}{}{4mm}  & \cr
          Q_c(k)/R && \to && H^1(k,F_{Q_c}) &\cr}$$ o\`u $F_{Z_c}$, resp. $F_{Q_c}$, est le
$k$-tore dual de
$\pic({\overline Z}_c)$, resp. $\pic({\overline Q}_c)$, et o\`u les fl\`eches horizontales sont
donn\'ees par les torseurs universels triviaux aux points marqu\'es \'evidents de $Z_c$ et
$Q_c$.  Sous les hypoth\`eses du th\'eor\`eme,  nous voulons montrer que la fl\`eche 
$Z_c(k)/R \to H^1(k,F_{Z_c})$ est surjective.
Par la th\'eorie des $k$-tores (Prop. 2.1 (iii)) on sait que la fl\`eche $ Q_c(k)/R  \to
H^1(k,F_{Q_c})$ est un isomorphisme et que la fl\`eche $Q(k) \to Q_c(k)/R$ est surjective.
On a vu dans la d\'emonstration du th\'eor\`eme 5.1 que $f^* : \pic({ \overline Q}_c) \to \pic({\overline Z}_c)$
est une injection scind\'ee de $g$-modules, dont le conoyau est un $g$-module stablement de
permutation. Le th\'eor\`eme 90 de Hilbert implique alors que
 la fl\`eche 
 $H^1(k, F_{Z_c}) \to H^1(k, F_{Q_c})$ est une bijection.
Pour $k$ un {\it  bon corps de dimension cohomologique $\leq 2$},
 la remarque 4.1.1
assure que la fl\`eche
 $Z(k) \to Q(k)$ est surjective. Ainsi pour tout tel corps la fl\`eche 
$Z_c(k)/R \to Q_c(k)/R \simeq H^1(k,F_{Q_c})$ est surjective. 
Ceci ach\`eve la d\'emonstration. \cqfd

\medskip

{\bf Question 6.2} {\it Sous les hypoth\`eses du th\'eor\`eme 6.1, la fl\`eche $ X_c(k)/R \to
H^1(k,F)$ est-elle une bijection ? }

\medskip

Le th\'eor\`eme 3.4 de [CTGiPa] montre que  ce
r\'esultat impliquerait
la finitude  de $ X_c(k)/R$ 
dans chacun des cas cit\'es \`a la remarque 4.1.1. La finitude 
dans le cas $p$-adique est facile \`a \'etablir, mais dans ce cas on aurait 
plus, on aurait la valeur exacte de $X_c(k)/R$.

La question 6.2 a une r\'eponse affirmative lorsque le sous-groupe ferm\'e connexe
$H$ est normal dans $G$ (Gille,  [G] et appendice de [B/K2]; Thm. 6.2 de [CT]).

Le premier cas \`a \'etudier est celui o\`u $k$ est un corps $p$-adique et $H$  un groupe
semi-simple. Dans ce cas, au vu de la proposition 3.1, la question est : 

{\it Soient $k$ un corps $p$-adique, $G$ un $k$-groupe semi-simple simplement connexe, $H
\subset G$ un $k$-sous-groupe ferm\'e connexe, semi-simple. Soit $X_c$ une $k$-compactification
lisse de $X=G/H$. La $R$-\'equivalence sur $X_c(k)$ est-elle triviale ?}

D'apr\`es un r\'esultat g\'en\'eral de Koll\'ar [K], la question
\'equivaut  \`a celle de la trivialit\'e de la $R$-\'equivalence sur l'ouvert $G/H = X \subset
X_c$.

\vskip1cm

{BIBLIOGRAPHIE}

\medskip

[Bog] F. A. Bogomolov, Groupe de Brauer des corps d'invariants de groupes alg\'ebriques (en russe), Mat. Sb.
{\bf 180} (1989) 279--293;  trad. ang.  Math. USSR-Sb. {\bf 66} (1990) 285--299.

[Bo1] M. Borovoi,  Abelianization of the second nonabelian Galois cohomology, Duke Math. J. {\bf
72} (1993) 217--239.

[Bo2] M. Borovoi, The Brauer-Manin obstructions for homogeneous spaces with connected or abelian
stabilizer, J. reine angew. Math. (Crelle) {\bf 473} (1996) 181--194.

[B/K1] M. Borovoi et B. Kunyavski\u{\i}, Formulas for the unramified Brauer group of a principal
homogeneous space of a linear algebraic group, J. Algebra {\bf 225} (2000) 804--821.

[B/K2] M. Borovoi et B. Kunyavski\u{\i}, Arithmetical birational invariants of linear algebraic
groups over two-dimensional geometric fields (with an appendix by P. Gille), J. Algebra  {\bf
276}  (2004) 292--339.

[CE] H. Cartan et S. Eilenberg, {\it Homological algebra}, Princeton, 1956.

[CT] J.-L. Colliot-Th\'el\`ene, R\'esolutions flasques des groupes r\'eductifs connexes, C. R.
Acad. Sc. Paris S\'er. I {\bf 339} (2004) 331--334.

[CTGiPa] J.-L. Colliot-Th\'el\`ene, P. Gille et R. Parimala, Arithmetic of linear algebraic
groups over 2-dimensional geometric fields, Duke Math. J. {\bf 121} (2004) 285--341.

[CT/K] J.-L. Colliot-Th\'el\`ene et B. Kunyavski\u{\i}, Groupe de Brauer non ramifi\'e des
espaces principaux homog\`enes des groupes lin\'eaires, J. Ramanujan Math. Soc. {\bf 13} (1998)
37--49.

[CT/San1] J.-L. Colliot-Th\'el\`ene et   J.-J. Sansuc,  La
$R$-\'equi\-valence sur les tores, Ann. Sc. E.N.S.  {\bf 10} (1977) 175--229.

[CT/San2] J.-L. Colliot-Th\'el\`ene et   J.-J. Sansuc, La descente sur les vari\'et\'es
rationnelles, II, Duke Math. J. {\bf 54} (1987) 375--492.

[CT/San3] J.-L. Colliot-Th\'el\`ene et J.-J. Sansuc, Principal homogeneous spaces under flasque
tori: applications, J.  Algebra {\bf 106} (1987) 148--205.

[CT/San4] J.-L. Colliot-Th\'el\`ene et  J.-J. Sansuc, The rationality problem for fields of
invariants under linear algebraic groups (with special regards to the Brauer group),
in {\it Proceedings of the Mumbai 2004 International Conference}, \`a para\^{\i}tre.

[D] A. Ducros, Dimension cohomologique et points rationnels,
J. Algebra {\bf 203} (1998) 349--354.

 [G] P. Gille, Cohomologie des groupes quasid\'eploy\'es sur des corps de dimension $\leq 2$,
Compositio Math. {\bf 125} (2001) 283--325.

[Gr] A. Grothendieck,  Le groupe de Brauer, I, II, III, in {\it Dix Expos\'es sur la Cohomologie
des Sch\'emas}, North-Holland, Amsterdam, 1968, pp. 46--188.

[K] J. Koll\'ar, Specialization of zero cycles, Publ. Res. Inst. Math. Sci.  {\bf 40}  (2004) 689--708.

[San] J.-J. Sansuc, Groupe de Brauer et arithm\'etique des groupes alg\'ebriques lin\'eaires sur
un corps de nombres, J. reine angew. Math. {\bf 327} (1981) 12--80.

[S] J-P. Serre, {\it Corps locaux}, deuxi\`eme \'edition, Hermann, Paris, 1968.

[Vos1] V. E. Voskresenski\u{\i},  Invariants birationnels des tores alg\'ebriques (en russe),  Uspehi Mat. Nauk  {\bf 30}  (1975), no. 2 (182), 207--208.

[Vos2] V. E. Voskresenski\u\i,  {\cyr Algebraicheskie tory} ({\it Tores alg\'ebriques}),  Nauka,  Moscou, 1977.  
  
[Vos3] V. E. Voskresenski\u{\i}, {\it Algebraic groups and their birational invariants},
Translations of Mathematical Monographs {\bf 179} (1998), American Mathematical Society.

\medskip

{\sevenrm

Jean-Louis Colliot-Th\'el\`ene,

C.N.R.S.,

UMR 8628, Math\'ematiques,
B\^atiment 425,
Universit\'e de Paris-Sud,

F-91405 Orsay

France

courriel : colliot@math.u-psud.fr

\bigskip

Boris  \`E.  Kunyavski\u{\i}

 Bar-Ilan University, Department of Mathematics 

52900 Ramat-Gan

Israel

e-mail : kunyav@macs.biu.ac.il }
\bye